%% file: arxiv_version.tex
\RequirePackage{tikz}
\documentclass[ijoc, nonblindrev]{arxiv_informs}
\OneAndAHalfSpacedXII

\usepackage{natbib}
\bibpunct[, ]{(}{)}{,}{a}{}{,}%
\def\bibfont{\small}%

\usepackage{lineno,hyperref}
\modulolinenumbers[5]
\usepackage{xcolor}
\newcommand{\william}[1]{\textcolor{blue}{#1}}

\usepackage[font=scriptsize]{subcaption}
\usepackage{booktabs}
\usepackage{amsmath, amsfonts, amssymb} 
\usetikzlibrary{arrows}
\usepackage{algpseudocode, algorithm}

\newtheorem{deff}{Definition}

\newcommand{\teal}[1]{\textcolor{teal}{#1}}
\newcommand{\margarida}[1]{\teal{ #1}}
\newcommand{\ie}{\emph{i.e.}}
\newcommand{\eg}{\emph{e.g.}}

\TheoremsNumberedThrough
\EquationsNumberedThrough 

\begin{document}


\TITLE{Adaptation, Comparison and Practical Implementation of Fairness Schemes in Kidney Exchange Programs}
\RUNTITLE{Adaptation, Comparison and Practical Implementation of Fairness Schemes in KEPs}





\RUNAUTHOR{St-Arnaud, Carvalho and Farnadi}

\ARTICLEAUTHORS{
\AUTHOR{William St-Arnaud}
\AFF{D\'epartement d'Informatique et de Recherche Op\'erationnelle \& CIRRELT, Universit\'e de Montr\'eal, 2920, chemin de la Tour, Montr\'eal, Quebec, H3T 1J4, Canada}
\AFF{Mila-Quebec AI Institute, 6666 St-Urbain, Montr\'eal, Quebec,  H2S 3H1, Canada}

\AUTHOR{Margarida Carvalho}
\AFF{D\'epartement d'Informatique et de Recherche Op\'erationnelle \& CIRRELT, Universit\'e de Montr\'eal, 2920, chemin de la Tour, Montr\'eal, Quebec, H3T 1J4, Canada}

\AUTHOR{Golnoosh Farnadi}
\AFF{Mila-Quebec AI Institute, 6666 St-Urbain, Montr\'eal, Quebec,  H2S 3H1, Canada}
\AFF{D\'epartement de sciences de la d\'ecision, HEC Montr\'eal, 3000, chemin de la Côte-Sainte-Catherine, Montr\'eal,  Quebec, H3T 2A7, Canada}
}

\ABSTRACT{
In Kidney Exchange Programs (KEPs), each participating patient is registered together with an incompatible donor. Donors without an incompatible patient can also register. Then, KEPs typically maximize overall patient benefit through donor exchanges. This aggregation of benefits calls into question potential individual patient disparities in terms of access to transplantation in KEPs. Considering solely this utilitarian objective may become an issue in the case where multiple exchange plans are optimal or near-optimal. In fact, current KEP policies are all-or-nothing, meaning that only one exchange plan is determined. Each patient is either selected or not as part of that unique solution. In this work, we seek instead to find a policy that contemplates the probability of patients of being in a solution. To guide the determination of our policy, we adapt popular fairness schemes to KEPs to balance the usual approach of maximizing the utilitarian objective. Different combinations of fairness and utilitarian objectives are modelled as conic programs with an exponential number of variables. We propose a column generation approach to solve them effectively in practice. Finally, we make an extensive comparison of the different schemes in terms of the balance of utility and fairness score, and validate the scalability of our methodology for benchmark instances from the literature.
}




\maketitle

\section{Introduction}
Kidney Exchange Programs (KEPs) allow patients that have an incompatible donor to exchange their donor’s kidney with the donor of another patient, thereby providing an additional alternative in the search for an organ. KEPs are implemented in several countries, \eg, South Korea \citep{park}, the Netherlands~\citep{deklerk}, the United Kingdom~\citep{manlove} and Canada~\citep{malik}.

The primary objective of KEPs is to find a solution, namely an exchange plan that consists of a set of donor exchanges resulting in compatible transplants for the patients. However, the main objective of maximizing the number of transplants does not consider the notion of fairness, which may result in some patients being excluded from the optimal exchange plan. To address this issue, several approaches to fairness in KEPs have been explored. These include maximizing the number of hard-to-match patients in a selected exchange plan \citep{dickerson2014balance} and balancing participation from different countries involved in a KEP pool \citep{shapley_kep, xenia2019, benedek2021computing, benedek2023partitioned}. These approaches are commonly referred to as \emph{group fairness}. Other fairness approaches involving individual patient probabilities of being selected in an exchange plan have been explored for pairwise matchings \citep{Roth_Sonme_Unver_2005_b, Jian2014} and more generally, for any maximum cycle length \citep{farnadi_if_kep}.

Of the latter type, one such approach is \emph{individual fairness} by \citet{farnadi_if_kep}, which uses probability distributions over exchange plans. These distributions can be viewed as lottery policies over the set of exchange plans. Lottery policies provide a more general framework to define fairness metrics than single (deterministic) exchange plans. A fair lottery policy has also been shown to be equivalent to a fair outcome \citep{bolton}. Maximizing the total number of transplants and optimizing a fairness score through deterministic policies may be at odds \citep{mcelfresh2018scalable}; therefore, it is crucial to consider lottery policies over all feasible exchange plans. For example, selecting a particular patient may severly impact the total number of transplants. It is still desirable not to fully exclude this patient from consideration. Choosing from the set of policies over all feasible exchange plans requires that the selected policy balances both objectives well, \ie, maximization of the expected number of transplants and maximization of a fairness score. This example motivates the need for a detailed study of fairness in KEPs to guide the computation of lottery policies. Hence, we formalize it under the lens of fairness schemes \citep{bertsimas} and we adapt (\ie, model) popular fairness scores to the context of KEP: individual fairness \citep{farnadi_if_kep}, Nash's standard of comparison \citep{nash, bertsimas}, Aristotle's principle of fairness \citep{aristotle}, and Rawlsian justice \cite{rawls_1973}. 
In this article, we aim to compare fairness schemes driven by a single-objective (a fairness score) and fairness schemes based on multi-objective models, which we refer to as the Social Welfare Program \citep{saaty1955} and the Nash Social Welfare Program \citep{Charkhgarda2020TheMO}. Since for each fairness scheme we build a lottery, our decision variables correspond to the probability of selecting an exchange plan. Thus, we first show that all our fairness schemes can be modeled as conic programs with an exponential number of variables and second, we devise a column generation to solve these programs effectively. As part of our experiments, we explore the quality of solutions returned by each framework-scheme combination in terms of the utilitarian (\ie number of transplants) and fairness objectives, as well as the effectiveness of our column generation methodology for their applicability to realistically sized KEP instances.
We can summarize our key contributions in this article as follows:
\begin{enumerate}
    \item The interpretation and modeling of popular fairness concepts (scores) in the context of KEPs. This exercise allows a unified view of fairness in KEPs, including the broadly-used group fairness.
    \item The construction of fairness schemes through the determination of lottery policies for KEPs based on a single-fairness score, and on the combination of utilitarian and fairness scores (bi-objective).
    \item For the bi-objective schemes, the use of the Social Welfare Program and the Nash Social Welfare Program to find a Pareto-optimal lottery policy.
    \item The formulation of our fairness schemes as conic programs and the design of a column-generation approach to solve them efficiently.
    \item A detailed experimental study comparing our fairness schemes, and providing insights to policy-makers on their benefits and limitations.
\end{enumerate}

\section{Related work}
\label{sec:related_work}


Fairness issues arise in several decision-making problems involving resource allocation. Common examples include fairness in communication networks, facility location and job scheduling \citep{hamoud}. For example, in the facility location problem, when minimizing the distance from facilities to service receivers, some facilities may end up far from some receivers and closer to others. This can be perceived to be unfair if we are speaking of facility types like schools or clinics. When devising a fairness criterion for a decision-making problem, one needs to ensure that it is chosen properly, as optimizing a fairness objective may even have the opposite intended effect and cause harm.

In the context of KEPs, the notion of \emph{group fairness} has been explored. It can take the shape of ensuring fairness towards pre-determined groups of patients, \eg, those that are deemed hard-to-match (\ie, highly-sensitized). However, this may negatively affect the total number of transplants realized under an exchange plan~\citep{dickersonPOF}. In this setting, the definition of the \textit{price of fairness} provides a measure of the cost of prioritizing these pre-determined groups to the total number of transplants~\citep{bertsimas, dickersonPOF}. \citet{dickersonPOF}~show that the price of fairness can become quite high when maximizing the number of highly-sensitized patients. They introduce a threshold to balance these two objectives: the threshold corresponds to the percentage of highly-sensitized patients that are to be included in the solution; it is given as a constraint and chosen empirically. \citet{freedman2018}~instead rely on maximizing the number of transplants and breaking ties using the number of highly-sensitized patients. The aforementioned fairness approaches in KEPs are not exhaustive. For examples related to game theory,  \citet{sonmez_Unver_2011, Ashlagi2014,Carvalho2019} and~\citet{xenia2019}~all cast the fairness component of KEPs in a framework where agents (representing hospitals, regions or countries) aim to maximize the benefit for their own patients; roughly, the goal is to find a policy encouraging agents to participate in a joint program. The various approaches to fairness in decision-making problems can be better formalized through the language of \textit{fairness schemes} \citep{bertsimas}.  We formalize fairness as an objective that is optimized. However, we also devise fairness schemes simultaneously optimizing a fairness objective and the usual utilitarian objective for KEPs in an attempt to balance the aggregated and individual patient benefits.

Also relevant within our context is literature related to \emph{procedural fairness}. \citet{bolton}~compare two fairness approaches that arise naturally in the context of games: fair procedures and fair outcomes. The authors argue that a fair procedure can be used as a substitute for a fair outcome. If we instantiate their definitions in KEPs, a procedure (policy) for the selection of an exchange plan would be considered fair if each patient considers their individual expected benefit satisfying. In such a procedure, a patient accepts the potential realization of an exchange plan (outcome) where their benefit is unsatisfying because in expectation, that is unlikely to occur. \citet{tyler_1996} also provides supporting evidence for this fact by demonstrating that a fair procedure is the determinant factor in the analysis of the fairness of an outcome. The study also shows that individuals' judgements about the fairness of a procedure is not influenced by the outcome if it was judged to be fair beforehand. The \emph{individual fairness} approach for KEPs found in~\citet{farnadi_if_kep} follows this line of research on fair procedures, although restricted to policies over optimal exchange plans, \ie, exchange plans achieving the maximum number of transplants. Their goal is to balance individual patient probabilities of being selected in an exchange plan. By definition, individual fairness entails a null \emph{price of fairness} to the detriment of maybe excluding some patients from consideration, \ie, those that cannot be part of an optimal exchange plan. However, one can also consider a fair procedure that consists of probability distributions over exchange plans that are not optimal in the (expected) number of transplants. This is what we set out to do in this article. Extending the support for our lotteries can allow us to give a non-zero probability of being in an exchange plan to patients who are not in optimal solutions. Moreover, we consider other popular fairness scores besides individual fairness, and we devise a methodology to solve real-sized instances; the methodology in \citet{farnadi_if_kep} is applicable to small KEP instances.

Finally, we will also consider fairness schemes combining the utilitarian objective with a fairness objective. This lead us to a bi-objective setting. A typical resolution approach is to seek a solution optimizing a linear combination of the objectives; this will be referred here as the Social Welfare Program (SWP)~\citep{saaty1955, Charkhgarda2020TheMO}. In \citet{Charkhgarda2020TheMO}, the authors introduce the framework of the Nash Social Welfare Program (NSWP) for multi-objective optimization. It addresses some weaknesses that are usually associated with the traditional Social Welfare Program approach, for instance, by being indifferent to the scale of the objectives. This framework is very similar to the previous work by \citet{wierzbicki_1980}, where the aim is to select a solution that is far from a reference point and that is also on the Pareto frontier; the key difference is that the latter method solves a min-max optimization problem by projecting a feasible point back on the Pareto frontier. We apply both the SWP and NSWP to the design of fairness schemes for KEPs, enabling us to evaluate experimentally their advantages when compared with single-objective fairness schemes.

\section{Kidney Exchange Programs}\label{sec:integerprog}
A KEP consists of a set of incompatible patient-donor pairs and a set of non-directed donors (NDD). Each incompatible patient-donor pair contains a patient necessitating a kidney and a (incompatible) donor willing to donate a kidney. This incompatibility can, for example, be due to bloodtype or tissue type mismatch \citep{canadian_blood_services_2019}. The program seeks to find compatible donor exchanges between its pairs and also with the NDDs. The objective is to form an exchange plan, in such a way that every patient relinquishing their donor also receives a donor in turn. This can be thought of as a form of bartering involving multiple parties.

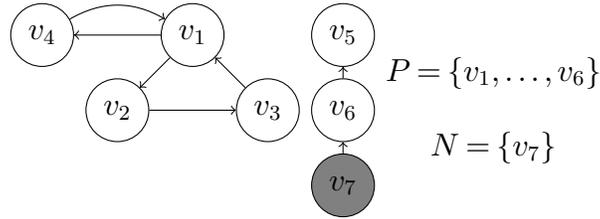
\begin{figure}[!ht]
\centering
\begin{tikzpicture}
\node[draw, circle] (v1) at (0,0) {$v_1$};
\node[draw, circle] (v2) at (-1,-1) {$v_2$};
\node[draw, circle] (v3) at (1, -1) {$v_3$};
\node[draw, circle] (v5) at (2,0) {$v_5$};
\node[draw, circle] (v6) at (2,-1) {$v_6$};
\node[draw, circle, fill=gray] (v7) at (2,-2) {$v_7$};
\node[draw, circle] (v4) at (-2,0) {$v_4$};

\draw[->] (v1) -- (v2);
\draw[->] (v2) -- (v3);
\draw[->] (v3) -- (v1);
\draw[->] (v7) to (v6);
\draw[->] (v6) to (v5);
\draw[->] (v4) to[bend left] (v1);
\draw[->] (v1) to (v4);

\node () at (4,-0.5) {$P = \{v_1,\dots,v_6\}$};
\node () at (4, -1.5) {$N = \{v_7\}$};
\end{tikzpicture}
\caption{Example of a KEP graph}
\label{fig:kep_ex}
\end{figure}

Formally, a KEP instance can be represented as a graph $G = (V,A)$, where $V = P \cup N$ consists of the set of incompatible patient-donor pairs $P$ and the set of non-directed donors $N$, and where $A$ is the set of arcs between these vertices. We have $(i,j) \in A$ if and only if the donor of vertex $i$ is compatible with the patient from vertex $j$. One can realize that no $(i,j) \in A$ exists for $j \in N$, \ie, no NDD can have an arc pointing to them. For an illustration, we refer the reader to Figure~\ref{fig:kep_ex}. Vertices $v_1$, $v_2$, $v_3$, $v_4$, $v_5$ and $v_6$ represent incompatible patient-donor pairs. They can swap donors to match their respective patients with a compatible organ. Vertex $v_7$ (shaded in gray) represents an NDD.

Exchange plans consist of arc-disjoint \emph{cycles} and \emph{chains}. A cycle is formed with incompatible patient-donor pairs, while a chain is a simple path starting in an NDD. The \emph{length} of a cycle (and chain) is the number of patients in it (\ie, number of vertices from P). In Figure~\ref{fig:kep_ex}, a feasible exchange plan can be given by the cycle $(v_1, v_2, v_3,v_1)$ of length 3 and the chain $(v_7, v_6, v_5)$ of length 2. An alternative exchange plan could instead involve the cycle $(v_4, v_1, v_4)$ and the same chain. For ethical reasons, the transplant operations corresponding to a selected cycle must be performed simultaneously. Therefore, due to transplant logistics, most KEPs impose a limit $K$ on the length of the cycles, and often, also a limit $K'$ on the length of the chains.

The typical goal of a KEP is the determination of an exchange plan maximizing the benefit of the patients. Thus, the simplest way of formulating this optimization mathematically is with integer programming through the \textit{cycle formulation}~\citep{Abraham2007, Roth2007}. 
For large instances, \textit{position-indexed edge formulations} \citep{dickersonPIEF} are generally more appropriate since we obtain a polynomial size to our mathematical model, while the cycle formulation can be exponentially large. A decomposition approach using branch-and-price has been explored by \citet{riascos} based on the work of \citet{constantino}. \citet{omer} provide a branch-and-price package to facilitate the deployment of related ideas while also contributing with an extensive analysis of the complexity of the pricing problems and suggesting algorithmic improvements to both the branch-and-price algorithms and the pricing problems.

Next, we describe a variant called \textit{hybrid position-indexed edge formulation}~\citep{dickersonPIEF} since it will be used within our decomposition methodology; the goal of identifying the best formulation is beyond the scope of this paper and we focus on the hybrid position-indexed edge formulation as it is both simple and not exponential in size. This formulation can be used directly with modern solvers and has been shown to solve realistic-sized instances. To simplify the exposition, we assume that $K'=K$. Here, for each $l \in V$, we create  a copy $G^l = (V^l, A^l)$ of the graph $G$. The sets $V^l$ contain all vertices with indices $i \geq l$ and the sets $A^l$ are the subsets of arcs induced by each $V^l$. We also define the following sets for $i,j, l \in V$:
\begin{align*}
    \mathcal{K}(i,j,l) &=
    \begin{cases}
    \{1\} &i=l \\
    \{2,\dots,K-1\} &i,j > l \\
    \{2,\dots,K\} &j=l.
    \end{cases} &
    \mathcal{K}'(i,j) &=
    \begin{cases}
    \{1\} &i \in N \\
    \{2,\dots,K\} &i \in P.
    \end{cases}
\end{align*}
The purpose of having many graph copies is to eliminate possible symmetries when expressing cycles. For graph copy $l$, we only have vertices with indices greater or equal to $l$. This amounts to selecting a unique starting point in a cycle (which will be vertex $l$ for the $l$-th copy). Because of the definition of $\mathcal{K}(i,j,l)$, we can recover the valid positions of each vertex in a cycle contained in graph copy $l$. The sets $\mathcal{K}'(i,j)$ serve the same function but correspond to chains.

We define the binary variables $x^l_{ijk}$ which take value 1 if the donation of a kidney from vertex $i$ to vertex $j$ appears in the $k$-th position of a cycle in the $l$-th graph copy. We also define the binary variables $y_{ijk}$ which take value 1 if the donation of a kidney from vertex $i$ to vertex $j$ appears in the $k$-th position of a chain. The formulation is as follows:
\begin{align*}    \tag{HPIEF}
    \label{eq:pief}
    &\max & \sum_{l \in P} \sum_{(i,j) \in A^l} \sum_{k \in \mathcal{K}(i,j,l)} w_{ij} x_{ijk}^l &+ \sum_{(i,j) \in A} \sum_{k \in \mathcal{K}'(i,j)} w_{ij} y_{ijk}\\
    &\text{s.t.} \\
    (c_1) & & \sum_{l \in P} \sum_{j : (j,i) \in A^l} \sum_{k \in \mathcal{K}(j,i,l)} x_{j,i,k + 1}^l &+ \sum_{j : (j,i) \in A^l} \sum_{k \in \mathcal{K}'(j,i)} y_{jik} \leq 1 & & \forall i \in V\\
    (c_2) & & \sum_{\substack{j : (j,i) \in A^l \, \land \\ k \in \mathcal{K}(j,i,l)}} x_{jik}^l &= \sum_{\substack{j : (i,j) \in A^l \, \land \\ k+1 \in \mathcal{K}(i,j,l)}} x_{i,j,k+1}^l &
    &\begin{cases}
        \forall l \in V \\
        \forall i \in \{l+1,\dots, \lvert V \rvert \} \\
        \forall k \in \{1,\dots, K-1\}
    \end{cases} \\
    (c_3) & & \sum_{j : (i ,j) \in A} y_{ij1} &\leq 1 & &\forall i \in N \\
    (c_4) & & \sum_{\substack{j : (j,i) \in A \, \land \, \\ k \in \mathcal{K}'(j,i)}} y_{jik} &\geq \sum_{j: (i,j) \in A} y_{i,j,k+1} & &
    \begin{cases}
        \forall i \in P \\ 
        \forall k \in \{1, \dots, K - 1\}
    \end{cases} \\
    & & x_{ijk}^l &\in \{0,1\} & 
    &\begin{cases}
        \forall l \in V \\
        \forall (i,j) \in A^l \\
        \forall k \in \mathcal{K}(i,j,l)
    \end{cases} \\
    & & y_{ijk} &\in \{0,1\} &
    &\begin{cases}
        (i,j) \in A \\
        k \in \mathcal{K}'(i,j)
    \end{cases}
\end{align*}
where $w_{ij}$ is the weight associated with the benefit to the patient in vertex $j$ when receiving a transplant from the donor in vertex $i$. Constraint $c_1$ in \eqref{eq:pief} ensures that no pair is part of more than one cycle or chain. Constraint $c_2$ (flow conservation) ensures that whichever pair gives a kidney also receives one in return. Constraint $c_3$ verifies that no NDD gives to more than one patient. Constraint $c_4$ allows a chain of current length $L < K$ to be extended by one at position $L+1$. In other words, if vertex $i$ gives to vertex $j$ at position $k+1$, vertex $i$ must necessarily have received an organ from some vertex at position $k$ in a chain. As a remark, since most KEPs' primary goal is to maximize the number of transplants, we generally have $w_{ij}$ equal to $1$. Non-unitary weights may correspond to utility values when receiving a kidney from a particular donor.


\section{Lottery and fair procedure: the role of distributions}
\label{sec:distributions}
With \ref{eq:pief}, a deterministic exchange is determined. Throughout this article, we take a broader perspective and we use probability distributions over exchange plans, \ie~lotteries. It might not be obvious at first how we can justify their use and even how this would work in practice. After all, only one exchange plan will be selected and performed. 

We first delineate between two notions: fair procedures and fair outcomes. A fair procedure in the context of resource allocation translates to giving a \textit{good} expected score to participants (\eg~in terms of utility, fairness measure, etc.). A fair outcome, rather, means that each participant gets their fair share with respect to a scoring mechanism. In KEPs, for example, we can think of the probability of receiving a transplant versus the outcome given by the selection of an exchange plan. For instance, a fair outcome would be an exchange plan maximizing the number of patients with O-blood type receiving  kidney. Meanwhile, guaranteeing this group of patients high odds of being part of an exchange plan corresponds to a fair procedure. It has been shown that a fair procedure can be used as a substitute for a fair outcome \citep{bolton}, \ie~individuals are indifferent in the choice between a fair procedure and a fair outcome. 

Based on this idea, we elect to focus on probability distributions over exchange plans, as individuals will be indifferent to this choice as long as their expected utility and fairness measures provided by the KEP are \emph{reasonable}. Also, this approach captures more complex notions of fairness that would not be applicable under a strictly deterministic setting. For example, we could look at each patient's probability of receiving a transplant and try to equalize these probabilities among patients. This would give better odds to certain patients that might systematically be left out under a utilitarian matching approach as they are hard to match. We also note that we will only consider patients that can be in at least one feasible exchange plan. This cannot be considered unfair as excluded vertices would never be selected in any scenario. We might need to adapt this idea in the dynamic case as these pairs could be matched with the arrival of future pairs. The dynamic case, however, is beyond the scope of this article and is only mentioned as an interesting research avenue, which can leverage in the work described here.

Based on~\citep{farnadi_if_kep}, we formally define a probability distribution over feasible exchanges. To this end, it is useful to think of a feasible exchange plan as a subgraph of the KEP graph. We call $\mathcal{F}_G$ the set of admissible subgraphs of $G$, \ie, subgraphs that correspond to feasible exchange plans. We will use the notation $V(S)$ for $S \in \mathcal{F}_G$ to mean the vertex set of subgraph $S$.
\begin{deff}
A probability distribution $\delta: \mathcal{F}_G \to \mathbb{R}_{\geq 0}$ over feasible exchanges for $G$ satisfies the condition $\sum_{S \in \mathcal{F}_G} \delta(S) = 1$. The set $\mathcal{F}_G(v) \subseteq \mathcal{F}_G$ denotes the set of feasible exchanges including pair $v \in P$ with  $\delta_v= \displaystyle \sum_{S \in \mathcal{F}_G(v)} \delta (S)$ being the probability given to $v$ by $\delta$.\label{def:probability_distribution}
\end{deff}
If we have such a probability distribution, we can draw a feasible exchange plan from it and implement this exchange plan in practice. Thus, pairs have a certain probability of being selected or not. Although the concept of utility for each transplantation $(i,j) \in A$ can mean various things (\eg, expected survival time after transplantation), in this article, we will consider the probability of receiving a transplant under a particular distribution over feasible exchanges. This is motivated by the fact that most programs privilege the number of transplants, \eg, see Figure 1 in~\citep{biro} compiling the objectives of KEPs in Europe. However, our methodology is applicable to any utilitarian function. Going back to our example in Figure~\ref{fig:kep_ex}, we see that there are multiple ways to select a cycle of a an exchange plan. We could consider the following two feasible exchange plans:
\begin{align*}
    S_1 &= (V, A_1)  & \textrm{where} \ \ & A_1= \{(v_1,v_2),(v_2,v_3),(v_3,v_1),(v_7,v_6),(v_6,v_5)\} \\
    S_2 &= (V, A_2) & \textrm{where} \ \ & A_2 = \{(v_4, v_1), (v_1, v_4),
    ,(v_7,v_6),(v_6,v_5)\}.
\end{align*}
We remark that $S_2$ does not contain a maximum number of vertices of $P$, but this would includes vertex $v_4$. However, it can make sense to consider a sub-optimal feasible exchange plan since it may represent a chance of being transplanted for a patient of a pair that would otherwise be excluded in the optimal solution (pair $v_4$ in Figure~\ref{fig:kep_ex}). A combination of the two solutions given by $\delta(S_1) = p, \delta(S_2) = 1 - p$ would give a probability $1-p$ to vertex $v_4$ and probability $p$ to vertices $v_2$ and $v_3$ of being selected in an exchange plan. Depending on the characteristics of the vertices that are selected, one solution could be considered fairer than the other. Furthermore, the process of improving $v_4$'s selection probability can be considered a fairness-enhancing approach.

\section{Fairness schemes}
\label{sec:fairness_schemes}
In this section, we describe the various fairness schemes that we will compare for KEPs, i.e., procedures to determine a lottery $\delta$. We do not claim to provide an exhaustive list of fairness schemes, but we will focus on the most well-known and foundational approaches~\citep{bertsimas}. Later, in Section~\ref{sec:fairnessuti}, we propose the application of the Social Welfare Program and the Nash Social Welfare Program to devise alternative fairness schemes for KEPs that combine those described in this section, in an attempt to mitigate their potential limitations.

\subsection{Mathematical description of fairness schemes}\label{subsec:backfairness}
In this section, we will discuss the concept of a fairness scheme \citep{bertsimas} in the context of KEPs.
\begin{deff}\label{def:utility_set}
Given an optimization problem of the form
\begin{align*}
    &\max \, \sum_{i=1}^m f_i(x, y) \\
    &\text{s.t.} \quad (x,y) \in X \subseteq \mathbb{R}^{n_1 + n_2},
\end{align*}
where $f: \mathbb{R}^{n_1 + n_2} \to \mathbb{R}^{m}$ and $\{1,\ldots, m\}$ is a set of agents, we define the \textbf{utility set} $U$:
\begin{equation*}
    U := \left\{ u \in \mathbb{R}^{m} \mid \exists (x,y) \in X \text{ s.t. } f(x, y) = u \right\}.
\end{equation*}
\end{deff}
In the context of KEPs, our set of agents are the patients and the optimization problem in Definition~\ref{def:utility_set} that we will use to build the utility set is going to be \eqref{eq:pief} with $$f_v(x, y)= \sum_{l \in V}\sum_{(i,v) \in A^l} \sum_{k \in \mathcal{K}(i,v,l)} w_{iv} x^l_{ivk} + \sum_{(i,v) \in A} \sum_{k \in K'(i,v)} w_{ij} y_{ivk} $$ for each $v \in P$ (here, $m = \lvert P \rvert$; we also slightly abused notation and index the function components with $P$). Again, variables $x$ and $y$ serve the same role as in \eqref{eq:pief}. By default, we use $w_{iv}=1$, \ie, $f_v(x, y)$ is 
$1$ if vertex $v$ is selected and $0$ otherwise. 
We will also need the definition of a \textit{fairness scheme}.
\begin{deff}
    A fairness scheme is any function $\mathcal{L}: 2^{\mathbb{R}^{m}} \to \mathbb{R}^{m}$ that takes as input a utility set $U$ for constraint set X and returns a utility vector. The notation $2^{\mathbb{R}^m}$ refers to the power set of $\mathbb{R}^m$, \ie, the set of subsets of $\mathbb{R}^m$.
\end{deff}

In practice, a fairness scheme often takes the form of some process or algorithm that selects a solution among the set of feasible solutions. This latter set is given implicitly through the input $U$, which describes the set of attainable utilities given the constraints $X$. In the typical setting where we maximize the number of transplants in a KEP exchange plan, we could see the function $\mathcal{L}$ as returning the utilities of an exchange plan maximizing said objective. Conceptually, this utilitarian objective is a fairness scheme. In spite of that, the goal of this work is to consider other fairness notions relevant to KEPs and afterward, derive a scheme balancing those notions and the utilitarian criterion.

When considering notions of fairness in a KEP, it is highly relevant to have some measure of the performance of the fairness scheme in terms of its cost to the maximum number of transplants, \ie, the cost to utility. For this, \citet{dickersonPOF, bertsimas} introduce the measure of \textit{price of fairness} (POF).
\begin{deff}
\label{def:pof}
    The price of fairness of a fairness scheme $\mathcal{L}$ for a utility set $U \subseteq \mathbb{R}^{m}$ is given by
    \begin{equation*}
        \frac{\sum_{i=1}^m u^*_i - \sum_{i=1}^m \mathcal{L}_i(U)}{\sum_{i=1}^m u^*_i},
    \end{equation*}
    where $u^* \in U$ is the maximum total utility that can be achieved as part of the set $U$, \ie, $\sum_{i=1}^m u^*_i \geq \sum_{i=1}^m u_i~\forall u \in U$. Note that $u^* = f(x^*, y^*)$ for $(x^*, y^*) \in X$, an optimal solution to the problem in Definition~\ref{def:utility_set}.
\end{deff}
A POF close to 0 means that the fairness scheme results in a low sacrifice of the maximum utility. Conversely, if it is 1, then we are sacrificing all the utility in order to achieve fairness. In Figure~\ref{fig:kep_ex}, we can see that any fairness scheme resulting from the selection of all the vertices in the graph, has POF equal to $0$. 
If we consider Figure~\ref{fig:kep_ex}, the suboptimal exchange plan $S_2$ leads to a POF equal to $0.2$.

\subsection{Utilitarian principle}
The \emph{utilitarian principle} \citep{young_equity:_1995, sen_economic_1997, bertsimas} is based on the idea that a fair allocation of resources maximizes social welfare, \ie, the aggregation of the agents' utility. The most common notion of social welfare in KEPs involves the number of transplants.\footnote{We could instead consider measures such us the expected number of years added to the patient's lifespan~\citep{krikov2007} or their associated quality adjusted life year~\citep{gloriePhd,GLORIE2021}.} 
As mentioned in Section~\ref{subsec:backfairness}, it is worth pointing out that the utilitarian principle can be recast in the language of fairness schemes as it meets the requirements of the definition. However, so as to avoid confusion, from now on, we will distinguish the utilitarian principle from the other fairness schemes by not referring to it as such, unless explicitly stated.

\subsection{Aristotle's equity principle}
\emph{Aristotle's equity principle} \citep{aristotle, bertsimas} is based on the idea that when competing for access to resources, a fair allocation should consider pre-existing rights to these resources. 
For example, a patient that has waited a long time in the KEP pool should receive preferential treatment in order to obtain a kidney. This prioritization also has the effect of minimizing the deterioration of patients' health \citep{CIHI2017}. Other patient features can also be weighted in order to define the pre-existing right. In KEPs, the panel-reactive antibodies scores (PRA) of patients are often considered as an indicator of priority. This score ranging between $0$ and $100$ measures the difficulty of matching a patient with a random donor from the population. A patient with a high PRA is called a hard-to-match.  The assumption is that the harder to match patients (pairs) will remain in the pool longer and thus may never be selected in time. Thus, the goal is to clear them from the pool as soon as there is a suitable match. In this scenario, the fairness scheme takes the form $\mathcal{L}^*(U) \in \underset{\mathcal{L}}{\arg\max} \sum_{v \in H} \mathcal{L}_v(U)$, where $H \subseteq P$ is the set of pairs with a PRA above a pre-determined threshold (\eg, PRA above 80\%) as used in~\citep{dickersonPOF}. This type of fairness notion is referred to as \textit{group fairness} in the literature~\citep{farnadi_if_kep}.

\subsection{Nash standard of comparison}
The \textit{Nash standard of comparison} \citep{nash, bertsimas} is based on the idea of $n$-player cooperative games. When $n=2$, it measures the relative change in each player's utility when transferring a small amount of resources to the opponent. If the relative gain to one player is greater than the relative loss to the other, then this transfer of resources is justified. The same concept can be generalized to more than two players.

\begin{deff}
Under the Nash standard of comparison (or proportional fairness), an allocation $\mathcal{L}(U)$ given by the fairness scheme $\mathcal{L}$ is fair if
\begin{equation}
    \sum_{i = 1}^m \frac{u_i - \mathcal{L}_i(U)}{\mathcal{L}_i(U)} \leq 0
    \label{eq:prop_fairness}
\end{equation}
for any $u \in U$.
\end{deff}
We measure the proportional change in utility when moving from $\mathcal{L}(U)$ to any other $u \in U$ and if the aggregate is always less than or equal to $0$, then the allocation $\mathcal{L}(U)$ is considered fair. If we consider a starting allocation vector for the patients vertices, finding a proportionally fair allocation roughly translates to a transfer of utility from high-utility vertices to low-utility vertices. For instance, suppose that  we start with a (deterministic) feasible exchange where some vertex $v$ is selected (\ie, $u_{v}=1$), and some vertex $v'$ is not selected (\ie, $u_{v'}=0$). Then, we determine a fairness scheme $\mathcal{L}(U)$, where only the utilities of $v$ and $v'$ are modified to $\mathcal{L}_{v}(U)=0.95$ and $\mathcal{L}_{v'}(U)=0.05$, \ie, the probability of $v$ and $v'$ to be selected is $95\%$ and $5\%$, respectively. This allocation $\mathcal{L}(U)$ would be justified. In fact, the $5\%$ increase represents a significant improvement for the low utility vertex, while a $95\%$ chance of transplantation is still remarkably high. It can be shown that the computation of a proportional fair scheme can be determined by maximizing the sum of the log-probabilities $\delta_v$ for the vertices (patients) belonging to some exchange plan. This is because the necessary and sufficient first-order optimality conditions for this solution correspond to the condition of the Nash standard of comparison (see \citep{bertsimas}). 

\subsection{Rawlsian justice}
\label{sec:fairnes_rawls}
In \citet{rawls_1973}, the author introduces the concept of the \emph{veil of ignorance}. The key assumption is that every individual is not made aware of their own situation, biases, view of the world, etc. The only aspect that remains is knowledge about the world and their rational mind. Rawls argues that under this veil of ignorance, individuals can negotiate the terms of a social contract that seems inherently fair.

One way to approach this idea is to maximize the utility of the least well-off agent. When dealing with distributions over feasible exchange plans, this amounts to maximize the minimum probability of being in a selected exchange plan among the vertices in $P$ of the KEP pool. For this to make sense, we can assume that $P$ does not contain vertices that cannot be part of any non-empty feasible exchange plan, \ie, they need to be matched in some $\mathcal{F}_G$. . Otherwise, the probability will always be 0, and we would get a skewed notion of fairness. 

As an example, we can refer the reader to Figure~\ref{fig:kep_ex}. Using the Rawlsian approach, to maximize the minimum probability over the vertices, we could select the cycle $(v_1,v_2,v_3,v_1)$ with probability $\frac{1}{2}$ and the cycle $(v_1,v_4, v_1)$ with probability $\frac{1}{2}$. Here, we consider probabilities over feasible exchanges to give a chance to vertex $v_4$. 


\subsection{Individual fairness}
\label{sec:IF}
$L_p$-fairness for KEPs was first introduced by \citet{farnadi_if_kep} as a measure of \emph{individual fairness}. It provides a notion of distance to the mean probability of being selected as part of a distribution over exchange plans. The lower the value of this metric, the more the patient probabilities of receiving a transplant are equalized. If all patients receive the same probability of receiving a transplant, the value will be 0.
\begin{deff}
The $L_p$-fairness of a probability distribution $\delta$ is given by 
\begin{equation}
    \sum_{v \in P} (\delta_v - \bar{\delta})^p, \label{obj:individual_fair}
\end{equation}
where $\bar{\delta} = \frac{1}{\lvert P \rvert} \sum_{v \in P} \delta_v$ is the mean probability of being selected as part of an exchange plan. A probability distribution $\delta$ is individually fair if it minimizes the Objective~\eqref{obj:individual_fair}.
\end{deff}

In the original setting, the support of $\delta$ is restricted to exchange plans maximizing the number of transplants. Individual fairness is seen to generalize the previous concept of \textit{Lorentz dominance} in pairwise KEPs, where the patient probabilities of being selected are equalized~\citep{Roth_Sonme_Unver_2005_b, Jian2014}. In this article, we do not restrict the support of our lottery $\delta$  as near-optimal exchange plans (in terms of utility) can be argued to be a good compromise with fairness. Moreover, we scale the use of individual fairness, as in \citet{farnadi_if_kep}, their approach only applied to small graph sizes.

\section{Approaches to MOO in KEPs}\label{sec:fairnessuti}
Previously, we underlined the importance of fairness in the selection of an exchange. However, in KEPs we still want to perform well in terms of utility. This motivates the present section, where we seek to to propose alternative fairness schemes combining the utilitarian scheme with each of the fairness schemes presented in the previous section. We must turn to multi-objective optimization to combine the two objectives. 

Multi-objective optimization (MOO) concerns problems of the form 
 \begin{align*}
     \max_{x\in X} \left( f_1(x), \ldots, f_k(x) \right),
     \tag{MOO}
     \label{Opt:multi}
 \end{align*}
 where $X \subseteq \mathbb{R}^n$ and $f_i$ for $i=1,\ldots,k$ are real-valued functions.  Typically, we are interested in the computation of \emph{Pareto-optimal} solutions. We call \emph{ideal} and \emph{nadir} the vectors corresponding to the best and worst value for each objective respectively (see formal definitions~\ref{def:ideal} and \ref{def:nadir} in the appendix).
 
\begin{deff}
A solution $x^\star \in X$ is Pareto-optimal to~\eqref{Opt:multi} if there is no $y \in X$ such that (i) $f_i(y) \geq f_i(x^\star)$ for all $i=1,\ldots,k$ and (ii) $f_j(y)>f_j(x^\star)$ for some $j \in \{1,\ldots,k\}$. The set of Pareto-optimal solutions is called Pareto frontier and the associated optimal values form the Pareto front.\label{def:pareto_opt}
\end{deff}

In Section~\ref{sec:swp}, we will review the most natural approach, called the Social Welfare Program, where we optimize a linear combination of objectives. In Section~\ref{sec:NSWP}, we move on to the Nash Social Welfare Program and discuss why it might be interesting compared to the previous approach. 

\subsection{The Social Welfare Program}
\label{sec:swp}
The Social Welfare Program (SWP)~\citep{saaty1955, Charkhgarda2020TheMO} provides a framework in the context of multi-objective optimization to analyze the trade-off between the objectives, which in our case are utility and fairness. Concretely, a Social Welfare Program for~\eqref{Opt:multi}  is any optimization problem of the form
\begin{equation*}
    \max_{x \in X} \sum_{i=1}^k \lambda_i f_i(x) 
\end{equation*}
where $\lambda_i > 0$ for $i=1,\ldots,k$ are fixed weights. 

From Definition~\ref{def:pareto_opt}, it is easy to see that a solution to an SWP is Pareto optimal for the respective multi-objective problem. Thus, by solving SWPs with different weights, various points of the Pareto frontier may be computed. In fact, since all our formulations in this article are convex, we can recover all points from the Pareto front using this approach. As discussed in~\citet{Charkhgarda2020TheMO}, there are inherent weaknesses that come with the SWP approach. The main weakness with which we shall be concerned in the context of this article deals with the non-existence of a criterion (or consensus) for selecting a Pareto-optimal solution. In particular, it is not obvious how to set the weights $\lambda$ to obtain a Pareto solution balancing our objectives. This weakness is a major one since we do not know how to compare two Pareto-optimal solutions. We refer the reader to \cite{Charkhgarda2020TheMO} for further details concerning weaknesses of the SWP and how they are addressed by the NSWP.

\subsection{The Nash Social Welfare Program}
\label{sec:NSWP}
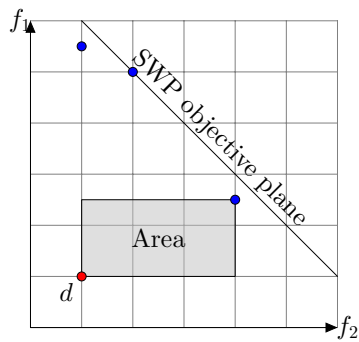
\begin{figure}[!ht]
\centering
\resizebox{.3\linewidth}{!}{
\begin{tikzpicture}[line cap=round,line join=round,>=triangle 45,x=1cm,y=1cm, line width=0.5pt]
\begin{scriptsize}
\draw[help lines] (-3,-3) grid (3,3);
\draw[->] (-3,-3) -- (3,-3);
\draw[->] (-3,-3) -- (-3,3);
\clip(-3.5,-3.5) rectangle (3.5,3.5);1

\draw[-] plot[domain=-2:3] (\x, {-\x + 1});

\draw[fill=lightgray, fill opacity=0.5] (-2,-2) rectangle (1, -0.5);
\draw[fill=blue] (-2,2.5) circle (2.5pt);
\draw[fill=blue] (-1,2) circle (2.5pt);
\draw[fill=blue] (1,-0.5) circle (2.5pt);

\draw[fill=red] (-2,-2) circle (2.5pt);

\draw (-2.3,-2.3) node () {\large $d$};
\draw (-3.2, 3) node () {\large $f_1$};
\draw (3.2, -3) node () {\large $f_2$};
\draw (-0.5, -1.25) node () {\large Area};
\draw (0.7, 0.7) node[rotate=-45] () {\large SWP objective plane};
\end{scriptsize}
\end{tikzpicture}
}
\caption[Pareto front $P^\star$ and the NSWP objective (gray area)]{Pareto front (blue dots) and the NSWP optimal objective (given by the gray area). The SWP optimal objective plane is depicted with $\lambda = (1, 1)$.\footnotemark}
\label{fig:pareto_nswp}
\end{figure}

The Nash Social Welfare Program (NSWP) \citep{Charkhgarda2020TheMO, rao1987} is an approach addressing the selection of a Pareto frontier point. It interprets the objective functions in Problem (MOO) as the payoffs of $k$ agents and it considers a reference $d_i \in \mathbb{R}, i=1,\ldots,k$ indicating to agent $i$ the minimum payoff they can guarantee. Hence, $f_i(x)-d_i$ becomes the surplus of agent $i$ for solution $x$, with the NSWP modelling the search for an $x \in X$ such that the multiplication of all agents' surplus is maximized. Mathematically, it has the form:
\begin{align*}\tag{NSWP}
\label{eq:NWSP}
\max &\prod_{i=1}^k (f_i(x) - d_i)^{\lambda_i} \\
\text{s.t.} \ \ & x \in X \\
&f_i(x) \geq d_i \quad i=1,\dots,k.
\end{align*}
The hyper-parameters $\lambda_i > 0$ are chosen beforehand and represent the weights attributed to each objective as in the SWP. In dimension $k=2$ with $\lambda = (1, 1)$, the NSWP attempts to maximize the total rectangular area defined by the reference point $d$ and the solution picked from
the Pareto front (see Figure~\ref{fig:pareto_nswp}). 

The NSWP possesses two key properties: it is both global-power-scale-free and local-benefit-scale-free. In other words, individually reweighting the terms $(f_i(x) - d_i)$ and/or globally rescaling the weights $\lambda_i$ will result in an equivalent problem up to some multiplication constant, \ie, the optimal solutions will be the same.
 \footnotetext{Inspired by figure in \citep{Charkhgarda2020TheMO}.}
We now proceed to the detailed modelling of the fairness schemes based on the NSWP, and the solution method to compute their respective lotteries effectively.

\section{Formulation and methodologies}\label{sec:ApplicationNSWP}
Let us stress the goals of the fairness schemes for KEPs to be developed in this section: (i) to balance between utility and a fairness concept,  (ii) to consider the chances of all patients to be in a selected exchange (fair procedure), and (iii) to efficiently run the scheme in practice. To achieve (i), we will apply SWP and~\ref{eq:NWSP}, and to attain (ii), we will define $X$ as the set of distributions over feasible exchanges. Since all distributions are given over the feasible exchanges, it is possible for a distribution to have an exponentially large support. This justifies goal (iii). Indeed, it might be that a distribution optimal to our SWP and~\ref{eq:NWSP} has a small support. Thus, we can instead seek to enumerate only the feasible exchange plans that will be part of the support of the optimal distribution. We will present a column generation method with this purpose. In this section, we lay out the conic programming formulations and column generation approach for the NSWP and  note that the fairness schemes of Section~\ref{sec:fairness_schemes} and the SWP are special cases of the NSWP.

For the sake of simplicity, in the following mathematical descriptions, we assume that the set $V$ only includes vertices that can be matched in at least one feasible exchange plan, \ie, $V = \cup_{S \in \mathcal{F}_G} V(S)$ (we abuse notation with the dual use of $V$ both as a set and a function). In practice, this set is computed before running the column generation: for each vertex $v$, we manually set Constraint $c_1$ of \eqref{eq:pief} to 1 (\ie, we turn it to an equality constraint) to force its inclusion in a solution. If one solution exists, $v$ is considered in our set $V$, otherwise, we can safely remove it from the graph along with incoming and outgoing arcs to this vertex. Assuming that $S \in \mathcal{F}_G$ is a feasible exchange plan, the variable $\delta_S$ will denote the probability of selecting this exchange plan, \ie, we use it as a short-hand for $\delta(S)$.

Next, we explicitly formulate SWP and~\ref{eq:NWSP} with $f_1$ given by the utilitarian principle and $f_2$ reflecting each of the following fairness concepts: individual fairness through the $L_1$-fairness, Aristotle's equity principle, Nash standard of comparison, and Rawlsian justice.\footnote{We begin with individual fairness since it uses most of the variables necessary to the other formulations.} We remark that we consider both objectives to have the same weight,  $\lambda_1=\lambda_2=1$, in~\eqref{eq:NWSP}. The idea behind this choice is that we wanted to give the same precedence to the utilitarian objective as the fairness objective. All mathematical formulations of this section will be convex programs.

\subsubsection*{Individual fairness} An individually fair scheme based on the $L_1$-fairness, minimizes this metric. Therefore, when in conjunction with the utilitarian principle,~\eqref{eq:NWSP} becomes the following conic quadratic program
\begin{align*}
    \label{eq:NSWP_l1}
    \tag{$\text{P}_{L_1}$}
    &&\max r \\
    \text{s.t.} \quad
    &(\alpha_0) & \sum_{S \in \mathcal{F}_G} \delta_S &= 1 \\
    &(\alpha_1) & y_1 &= \sum_{S \in \mathcal{F}_G} \sum_{v \in V(S)\setminus N} \delta_S - d_1\\
    &(\alpha_2) & y_2 &= -T - d_2 \\
    &(\alpha_3) & |P'| z &= \sum_{S \in \mathcal{F}_G} \sum_{v \in V(S)\setminus N} \delta_S \\
    &(\beta_v) & z_v &= \sum_{S \in \mathcal{F}_G \mid v \in V(S)\setminus N} \delta_S - z &\forall v \in P' \\
    &(\lambda_S) & \delta_S &\geq 0 &  \forall S \in \mathcal{F}_G \\
    &(u \in \mathcal{Q}_{r}^3) & (y_1,y_2,r) &\in \mathcal{Q}_{r}^3 \\
    &(w_v \in \mathcal{Q}^2) & (t_v,z_v) &\in \mathcal{Q}^2 &\forall v \in P' \\
    &(\eta) & \sum_{v \in P'} t_v &= T,
\end{align*}
where in parenthesis we denote dual variables and $P'=V \cap P$ is the set of pairs covered by some feasible exchange plan (remember the assumption on $V$ mentioned above). The computation of the reference  point $d$ is discussed in Appendix~\ref{sec:appendix_reference}. The objects $\mathcal{Q}^2$ and $\mathcal{Q}_{r}^3$ are the second-order and rotated second-order cones respectively (see Appendix~\ref{sec:appendix_def} for the definitions).

The constraint corresponding to dual variable $\alpha_0$ ensures that we have a probability distribution over feasible exchange. The variables $y_1$ and $y_2$ represent each multiplicative term in~\eqref{eq:NWSP}: $y_1$ is the expected number of transplants minus $d_1$ (constraint with dual variable $\alpha_1$) and $y_2$ is the negative of the $L_1$ metric minus $d_2$ since we are minimizing this loss (constraint with dual variable $\alpha_2$). The variable $z$ is equal to the the average pair probability of being selected in a an exchange plan (constraint with dual variable $\alpha_3$) and $z_v$ is the difference between $z$ and $\delta_v$ (constraint with dual variable $\beta_v$; recall Definition~\ref{def:probability_distribution}). Through $\mathcal{Q}^2$, we ensure that $t_v$ is equal to the absolute value of $z_v$ (constraint with dual variable $w_v$). Through $\mathcal{Q}_{r}^3$,  we guarantee that the optimal objective value $r$ coincides with $2 \cdot y_1 \cdot y_2$ which is equivalent to maximizing $y_1 \cdot y_2$ (constraint with dual variable $u$).

This problem formulation can have exponentially many variables $\delta_S$, so feasible exchanges $S$ will be generated as needed in a column generation fashion. To do so, we take the dual formulation of \eqref{eq:NSWP_l1}:
\begin{align*}
    \label{eq:dual_NSWP_l1}
    \tag{$\text{D}_{L_1}$}
    &&\min -\alpha_0 + \alpha_1 d_1 + \alpha_2 d_2 \\
    &\text{s.t.} 
    &\begin{pmatrix}
    \alpha_1 \\
    \alpha_2 \\
    1
    \end{pmatrix} + u &= 0 \\
    &&\begin{pmatrix}
    \eta \\
    \beta_v
    \end{pmatrix} + w_v &= 0 &\forall v \in P' \\
    &&\alpha_2 - \eta &= 0 \\
    &&\sum_{v \in P'} \beta_v + |P'|\alpha_3 &= 0 \\
    &&\alpha_0 - \sum_{v \in V(S) \setminus N} (\alpha_1 + \alpha_3 + \beta_v) + \lambda_S &= 0 &\forall S \in \mathcal{F}_G \\
    &&u &\in \mathcal{Q}_{r}^3 \\
    &&w_v &\in \mathcal{Q}^2 &\forall v \in P' \\
    &&\alpha \in \mathbb{R}^{4}, \beta \in \mathbb{R}^{|P'|}, \eta \in \mathbb{R}, \lambda & \geq 0.
\end{align*}
Using \eqref{eq:dual_NSWP_l1}, we immediately see from $\lambda_S \geq 0$ that we can rewrite the Constraints~\eqref{eq:NSWP_subprob_obj}:
\begin{align*}
    \label{eq:NSWP_subprob_obj}
    \tag{$\dag_S$}
    &-\alpha_0 + \sum_{v \in V(S)\setminus N} (\alpha_1 + \alpha_3 + \beta_v) \geq 0 &\forall S \in \mathcal{F}_G.
\end{align*}

In this way, given an optimal solution of~\eqref{eq:NSWP_l1} restricted to a subset of $\mathcal{F}_G$ (restricted master program), we can verify its optimality, or equivalently, dual feasibility, by finding the exchange plan minimizing the left-hand-side of Constraint~\eqref{eq:NSWP_subprob_obj}. Thus, in Algorithm~\ref{alg:columngen}, we use Constraint~\eqref{eq:NSWP_subprob_obj} as the objective of the subproblem during the solution generation process. Indeed, the subproblem implements this objective within~\eqref{eq:pief}. Our methodology stops, once the objective value of the subproblem is non-negative, \ie, no Constraint~\eqref{eq:NSWP_subprob_obj} is violated. In essence, Algorithm~\ref{alg:columngen} defines a cutting plane procedure for the dual problem~\eqref{eq:dual_NSWP_l1} and, consequently, a column generation for the primal~\eqref{eq:NSWP_l1}. Since the problem of interest is~\eqref{eq:NSWP_l1}, we refer to our methodology as a column generation.

In the remainder of this section, we do not provide details about the column generation approach for the next formulations given their similarity with the one just presented.

\begin{algorithm}
\begin{algorithmic}
\Function{MasterProblem}{$G$}
\State $\mathcal{F'}_G \gets \emptyset$ \Comment{Initialize the set of solutions}
\State $\zeta \gets -1$
\While{$\zeta <0$}
\State $\mathcal{D} \gets$~\eqref{eq:dual_NSWP_l1} restricted to Constraints~\eqref{eq:NSWP_subprob_obj} for $\mathcal{F'}_G$ 
\State $(\alpha^\star, \beta^\star, \lambda^\star, u^\star, v^\star, \eta^\star) \gets $ get solution of $\mathcal{D}$
\State $S, \zeta \gets$ {\sc Subproblem}$(\alpha^\star, \beta^\star, \lambda^\star, u^\star, v^\star, \eta^\star)$
\State $\mathcal{F'}_G  \gets \mathcal{F'}_G  \cup \{S\}$ \Comment{Store new solution/column}
\EndWhile
\State \Return $\mathcal{F'}_G$ and the associated primal variables $\delta^\star$
\EndFunction
\State
\Function{Subproblem}{$\alpha^\star, \beta^\star, \lambda^\star, u^\star, v^\star, \eta^\star$} \Comment{Solve the subproblem}
\State Set $w_{iv}=\alpha_1^\star+\alpha_3^\star+\beta_v^\star$ for $v \in P'$ in~\eqref{eq:pief}
\State Solve~\eqref{eq:pief} and denote by $S^\star$ the optimal solution
\State $V(S^*) \gets \{ v \in V \mid v \text{ is selected in the exchange plan } S^\star \}$
\State \Return $S^\star, \alpha_0^\star - \sum_{v \in V(S^\star)\setminus N} (\alpha_1^\star + \alpha_3^\star + \beta_v^\star)$
\EndFunction
\end{algorithmic}
\caption{Decomposition algorithm to solve~\eqref{eq:NSWP_l1}}
\label{alg:columngen}
\end{algorithm}

\subsubsection*{Aristotle's equity principle} One example of Aristotle's equity principle for KEPs involves the use of hard-to-match patients, \ie, patients with a PRA above a certain threshold. Here, we use $80\%$ as the threshold: every patient above this percentage is labelled as hard-to-match and in turn, its corresponding vertex (as done in~\citet{dickersonPOF}). We use the notation $V_H = P_H \cup N$, where $P_H = \{v \in P' \mid v \text{ has PRA} \geq 80\% \}$, to denote the restriction of the vertex set $V$ to pairs with high PRA (and in $P'$) and NDDs.
\begin{align*}
    \label{eq:NSWP_aristotle}
    \tag{$\text{P}_{\text{Aristotle}}$}
    &&\max r \\
    \text{s.t.} \quad
    &(\alpha_0) & \sum_{S \in \mathcal{F}_G} \delta_S &= 1 \\
    &(\alpha_1) & y_1 &= \sum_{S \in \mathcal{F}_G} \sum_{v \in V(S) \setminus N} \delta_S - d_1\\
    &(\alpha_2) & y_2 &= \sum_{S \in \mathcal{F}_G} \sum_{v \in V(S) \cap P_H} \delta_S - d_2 \\
    &(\lambda_S) & \delta_S &\geq 0 & \forall S \in \mathcal{F}_G\\
    &(u \in \mathcal{Q}_{r}^3) & (y_1,y_2,r) &\in \mathcal{Q}_{r}^3. 
\end{align*}
All variables keep analogous meanings to the ones introduced for~\eqref{eq:NSWP_l1}, except for $y_2$ which is now the expected number of hard-to-match patients in a selected exchange plan. The dual formulation \eqref{eq:dual_NSWP_aristotle} can be found in the appendix.

\subsubsection*{Nash standard of comparison} When using the Nash standard of comparison, we seek to maximize the sum of the log-probabilities of being selected in an exchange plan, \ie, $\sum_{v \in P'} \log \delta_v$. We then have the following~\eqref{eq:NWSP} formulated as an exponential conic program:
\begin{align*}
    \label{eq:NSWP_sumlog}
    \tag{$\text{P}_{\text{Nash}}$}
    &&\max r \\
    \text{s.t.} \quad 
    &(\alpha_0) & \sum_{S \in \mathcal{F}_G} \delta_S &= 1 \\
    &(\alpha_1) & y_1 &= \sum_{S \in \mathcal{F}_G} \sum_{v \in V(S) \setminus N}  \delta_S - d_1\\
    &(\alpha_2) & y_2 &= T - d_2 \\
    &(\beta_v) & z_v &= \sum_{S \in \mathcal{F}_G \mid v \in V(S)} \delta_S &\forall v \in P' \\
    &(\eta) & \sum_{v \in P'} t_v &= T 
    \\
    &(\lambda_S) & \delta_S &\geq 0 & \forall S \in \mathcal{F}_G \\
    &(w_v \in (K_{\text{exp}})^*) & (z_v,1,t_v) &\in K_{\text{exp}} &\forall v \in P'\\
    &(u \in \mathcal{Q}_{r}^3) & (y_1, y_2, r) &\in \mathcal{Q}_{r}^3,
\end{align*}
 where $K_{\text{exp}}$ is the exponential cone (defined in Appendix~\ref{sec:appendix_def}).

The variable $z_v$ now corresponds to the probability that $v$ is selected in an exchange plan. The cone $K_\text{exp}$ ensures that the variable $t_v$ is equal to $\log z_v$ in order to compute the sum of the logarithmic terms. The dual formulation \eqref{eq:dual_NSWP_nash} can be found in the appendix.

\subsubsection*{Rawlsian justice} Using the Rawlsian approach, we aim to balance the utilitarian principle with the minimization of the probability of a vertex to be selected in an exchange plan. The corresponding~\eqref{eq:NWSP} formulation is
\begin{align*}
    \label{eq:NSWP_mp}
    \tag{$\text{P}_{\text{Rawls}}$}
    &&\max r \\
    \text{s.t.} \quad
    &(\alpha_0) & \sum_{S \in \mathcal{F}_G} \delta_S &= 1 \\
    &(\alpha_1) & y_1 &= \sum_{S \in \mathcal{F}_G} \sum_{v \in V(S) \setminus N} \delta_S - d_1\\
    &(\alpha_2) & y_2 &= T - d_2 \\
    &(\beta_v) & z_v &= \sum_{S \in \mathcal{F}_G \mid v \in V(S) } \delta_S &\forall v \in P' \\
    &(\gamma_v) & z_v &\geq T &\forall v \in P' \\
    &(\lambda_S) & \delta_S &\geq 0 & \forall S \in \mathcal{F}_G\\
    &(u \in \mathcal{Q}_{r}^3) & (y_1,y_2,r) &\in \mathcal{Q}_{r}^3.
\end{align*}

Variable $T$ will be made equal to the minimum of variables $z_v$ (vertex probabilities of being selected in an exchangep plan) through the constraints identified with dual variables $\gamma_v$. The dual formulation \eqref{eq:dual_NSWP_mp} can be found in the appendix. If we are only interested in the fairness schemes of Section~\ref{sec:fairness_schemes}, we can simply replace the constraints with the dual variable $\alpha_1$ to $y_1 = 1$ in all the formulations provided in this section. If we are interested in the SWP, we can replace $r$ by $\lambda_1 * y_1 + \lambda_2 * y_2$ and remove the constraints associated with the variable $u$.

\section{Evaluation}\label{sec:experiments}

Next, we analyze multiple experiments designed to evaluate the benefits of introducing fairness schemes from Section~\ref{sec:fairness_schemes} and the fairness schemes based on the SWP and the NSWP. Prior to this, we describe our computational environment, instances and tested schemes.


\subsection{Experimental setup}

Each experiment was run on one cluster node powered by two Intel E5-2683 v4 Broadwell CPUs running at 2.1Ghz. A time limit of 1 hour per experiment was set, and 16GB of memory was used. All the code was implemented in Julia 1.7.1 and optimization problems were solved with Mosek 9.2.\footnote{\url{www.mosek.com}} A tolerance parameter of $10^{-8}$ was used when solving the NSWP in Mosek.

Based on the Saidman generator~\citep{saidman}, \citet{DickersonInstances} created graph instances (including PRA) that are available in Preflib~\citep{preflib2013}. 
Table~\ref{tab:instances}  summarizes the sizes and number of instances used by us from this repository. These sizes range over some KEP pool sizes used in Europe~\citep{biro_building_2018} and Canada~\citep{KPD_report2018}. 
The cap limit for cycles and chains length was set to $K=3$ and $K'=3$. These limits are consistent with some existing KEPs \citep{biro_building_2018}, and increasing them would mainly augment the set of feasible exchange plans available for our fairness schemes. 

\begin{table}[!ht]
    \centering
    \resizebox{0.25\linewidth}{!}{
    \begin{tabular}{c|c|c|c|c}
        \toprule
         &  \multicolumn{4}{c}{Percentage of NDDs}  \\
         $\lvert P \rvert$ & 0\% & 5\% & 10\% & 15 \% \\
         \midrule
         16 & 10 & 0 & 10 & 10\\
         32 & 10 & 10 & 10 & 10 \\
         64 & 10 & 10 & 10 & 10 \\
         128 & 10 & 10 & 10 & 10 \\
         256 & 10 & 10 & 10 & 10 \\
         \bottomrule
    \end{tabular}}
    \caption{Number of instances for each combination of $\lvert P \rvert$ and percentage of NDDs}
    \label{tab:instances}
\end{table}
We did not consider larger instances as we either did not have enough computing power or memory to build the~\eqref{eq:pief} formulation. Note that the efficiency of our decomposition method goes hand in hand with the state-of-the-art (integer programming) approaches for KEPs. These correspond to the subproblems in our column generation procedure (recall Algorithm~\ref{alg:columngen}), whose improvement is beyond the scope of this article.

\subsection{Results}

In the experimental results, we abbreviate to Utilitarian, Nash, Rawls and IF, the fairness schemes considering the utilitarian principle (number of transplants), Nash standard of comparison, Rawlsian justice and individual fairness, respectively. We use the term ``Single'' when a scheme is based on the optimization of a single fairness concept, ``NSWP'' when the NSWP scheme is considered and ``SWP'' when the objective is a linear combination of a fairness concept and the utilitarian objective. The weights for the linear combination were chosen as $(\lambda_1,\lambda_2) = \left (\frac{1}{i_1 - d_1}, \frac{1}{i_2 - d_2} \right)$ in order to have similar scales for both objectives, with $i$ and $d$ being the ideal and reference point, respectively. See Appendix~\ref{sec:appendix_reference} for details on the computation of $d$.

The set of research questions that we investigate can be divided into two groups:
\begin{enumerate}
    \item Research questions 1, 2, and 3 analyze the effect of combining a utilitarian objective with a fairness objective. We analyze how effective the methods are in terms of the quality of the solutions with respect to each fairness objective.
    \item Research questions 4 and 5 empirically measure the efficiency of the column generation approach (Algorithm~\ref{alg:columngen} for NSWP) in practice.
\end{enumerate}

We provide complementary research questions in Appendix~\ref{sec:appendix_dist} and \ref{sec:appendix_pou}.

\noindent\textbf{RQ1: What is the effect of different fairness schemes on utility?} 
First, we compute the average POF of each fairness scheme over all instances. These values can be found in Table~\ref{tab:objective_cost}. Only instances for which all schemes were solved within the time limit were considered (details provided in Appendix~\ref{sec:appendix_solved}).

\begin{table}[!ht]
    \centering
    \resizebox{0.6\textwidth}{!}{
    \begin{tabular}{c|c|c|c|c}
        \toprule
        & \multicolumn{4}{c}{Fairness concept} \\
        Scheme & IF & Rawls & Aristotle & Nash\\
        \midrule
        Single & $0.408 \pm 0.461$ & $0.089 \pm 0.150$ & $0.074 \pm 0.190$ & $0.007 \pm 0.021$\\
        SWP & $0.100 \pm 0.147$ & $0.014 \pm 0.044$ & $0.000 \pm 0.000$ & $0.004 \pm 0.014$ \\
        NSWP  & $0.118 \pm 0.157$ & $0.005 \pm 0.015$ & $0.001 \pm 0.005$ & $0.004 \pm 0.014$ \\
        \bottomrule
    \end{tabular}}
    \caption{POF for each scheme (all graph sizes)}
    \label{tab:objective_cost}
\end{table}
 As expected, the Single schemes  have the highest POFs in comparison with SWP and NSWP since the number of transplants is not considered. This supports that schemes based solely on fairness may not be adequate, since they may significantly decrease the utilitarian value. Overall, the NSWP has on average very similar POFs  to that of the SWP schemes. The POFs are low indicating that fairness does not come at the cost of utility. IF is the only exception to the similarity between SWP and NSWP. However, we note that this is a byproduct of the fact that the empty (or a non-maximal) exchange plan may be selected in NSWP (see Appendix~\ref{sec:appendix_bilevel}).

\begin{figure}[!ht]
    \centering
    \begin{subfigure}{0.45\textwidth}
        \centering
        \includegraphics[width=\textwidth]{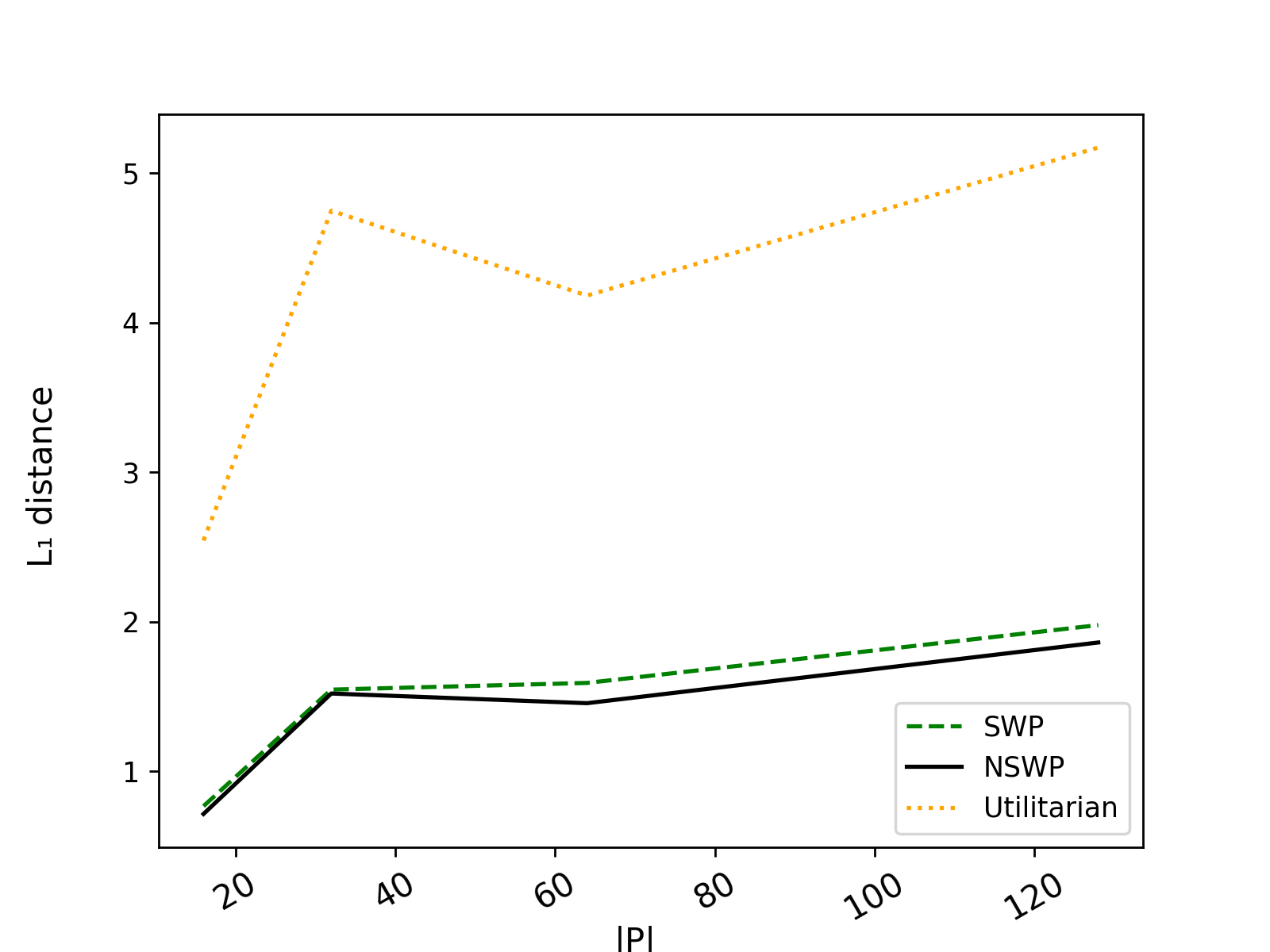}
        \caption{IF}
    \end{subfigure}
    \begin{subfigure}{0.45\textwidth}
        \centering
        \includegraphics[width=\textwidth]{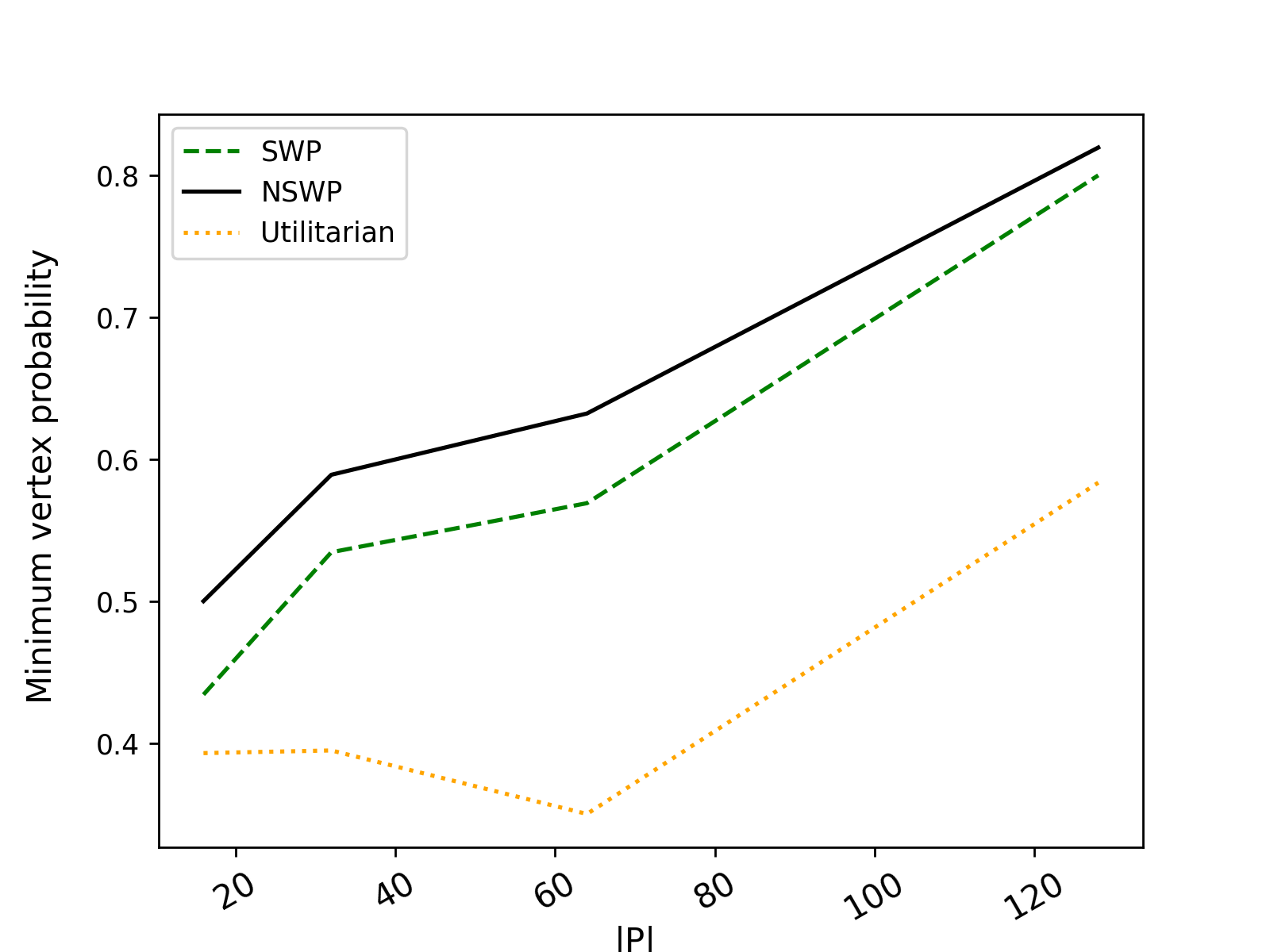}
        \caption{Rawls}
    \end{subfigure}

    \begin{subfigure}{0.45\textwidth}
        \centering
        \includegraphics[width=\textwidth]{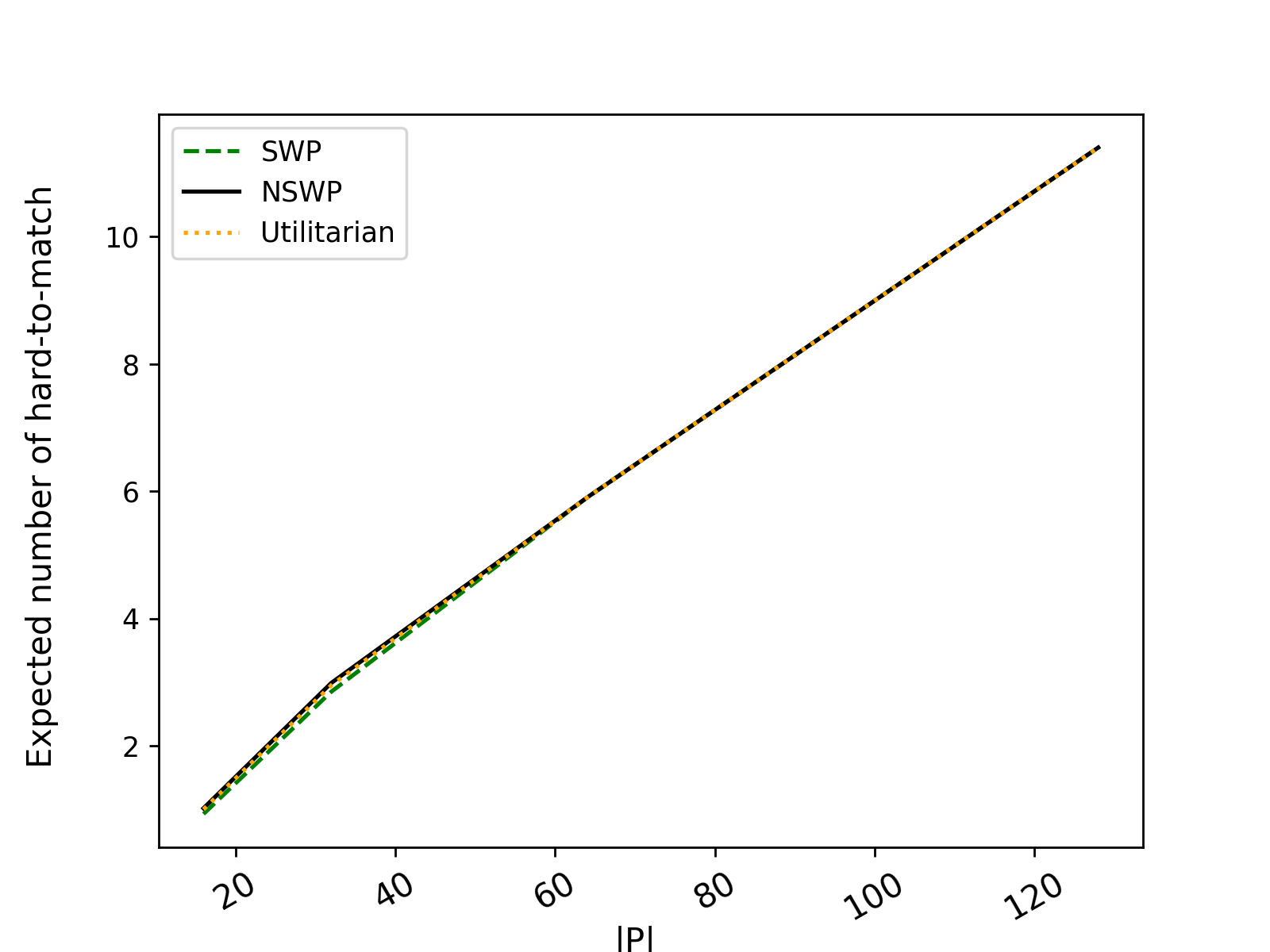}
        \caption{Aristotle}
    \end{subfigure}
    \begin{subfigure}{0.45\textwidth}
        \centering
        \includegraphics[width=0.98\textwidth]{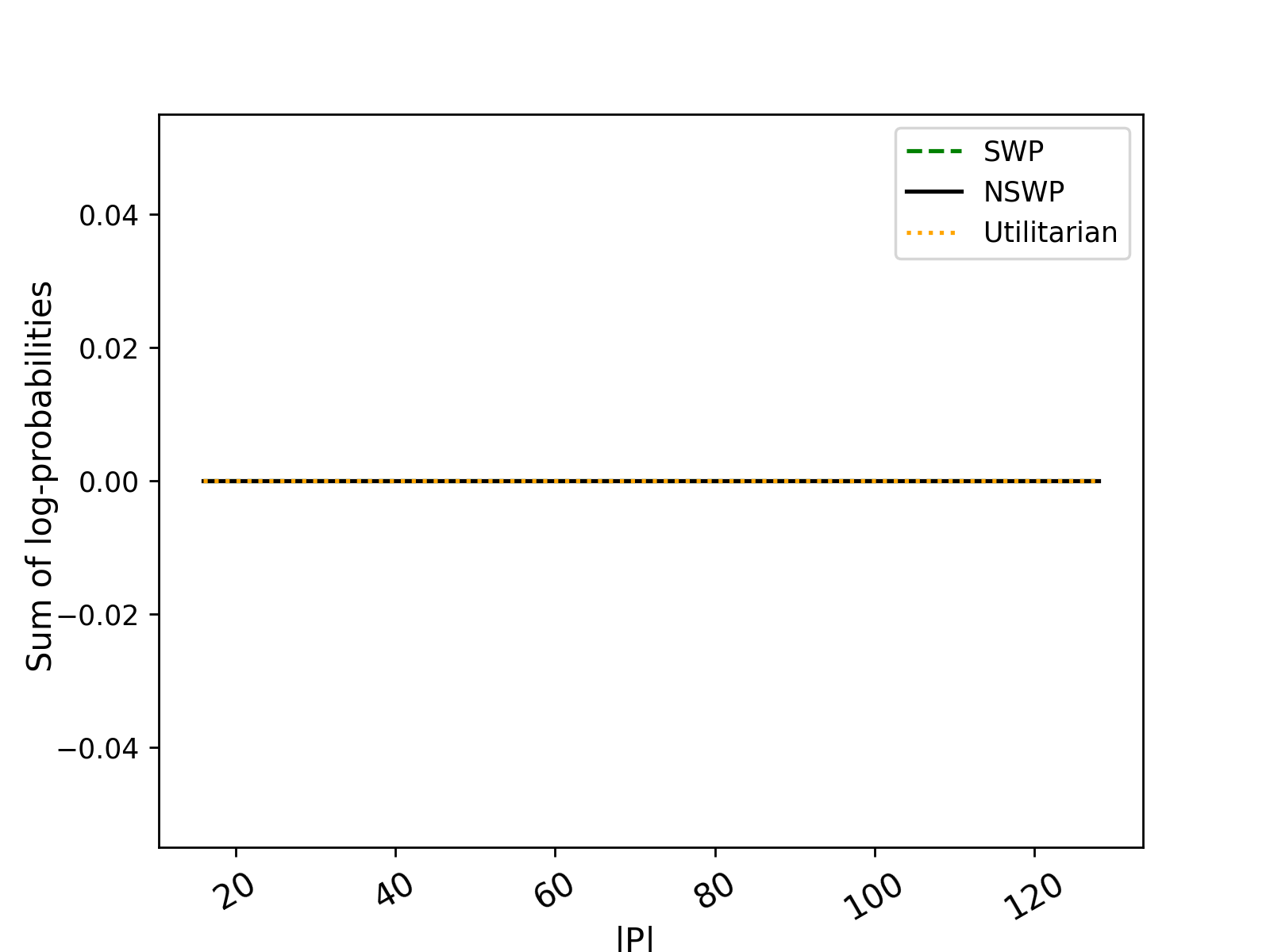}
        \caption{Nash}
    \end{subfigure}
    \caption{Comparison of SWP, NSWP and the utilitarian approach in terms of the value for each fairness objective}
    \label{fig:POU}
\end{figure}

\noindent\textbf{RQ2: What is the effect of the utilitarian approach on the fairness objectives?} 
Given the good POF of the SWP and the NSWP, in this question, we focus on those schemes. We looked at the value of the fairness objective for the NSWP, SWP and the utilitarian approach. The goal is to observe how much of an improvement is obtained by the NSWP or SWP in the fairness value when compared to the standard utilitarian method. The values can be found in Figure~\ref{fig:POU}. Again, we only consider instances for which we were able to apply the three approaches within the time limit.

A marked difference can be observed for IF and Rawls, where the NSWP and SWP show a significant improvement to the fairness objective. However, we do not see any difference for Aristotle and Nash's fairness principles. The Aristotelian case is made obvious by our analysis of the support size in RQ5. Essentially, for many instances there is an exchange plan achieving the ideal values. For Nash's fairness objective, we averaged the values that did not return negative infinity, a recurrent property for the utilitarian scheme as some vertices are not included in the support of the NSWP and SWP. 

\noindent\textbf{RQ3: What is the minimum number of pairs included in exchange plans from the support?} For this question, we investigate how acceptable the selection of an exchange plan would be to patients affected by our methodology. Given a lottery, \ie, a probability distribution, we look at the exchange plans with strictly positive probability. Among those, we compute the minimum number of pairs that are part of a single exchange plan. Table~\ref{tab:support_min} reports the values obtained as a fraction of the the maximum number of pairs that can be included as part of the utilitarian approach.

\begin{table}[!ht]
    \centering
    \resizebox{0.6\textwidth}{!}{
    \begin{tabular}{c|c|c|c|c}
        \toprule
        & \multicolumn{4}{c}{Fairness concept} \\
        Scheme & IF & Rawls & Aristotle & Nash\\
        \midrule
        Single & $0.489 \pm 0.491$ & $0.566 \pm 0.442$ & $0.905 \pm 0.205$ & $0.965 \pm 0.084$\\
        SWP & $0.551 \pm 0.44$ & $0.697 \pm 0.407$ & $1.000 \pm 0.005$ & $0.965 \pm 0.084$ \\
        NSWP  & $0.541 \pm 0.449$ & $0.806 \pm 0.261$ & $0.995 \pm 0.032$ & $0.967 \pm 0.083$ \\
        \bottomrule
    \end{tabular}}
    \caption{ Minimum fraction of patients among support for each scheme (all graph sizes)}
    \label{tab:support_min}
\end{table}

We can observe that the SWP and the NSWP perform similarly in terms of the minimal fraction of patients that are included among the solutions in the support (\ie, with strictly positive probability). Both provide a noticeable improvement over the single objective case for Rawls and Aristotle and slightly less so for IF. The main takeaway is that both SWP and NSWP always perform at least as well as the single objective case on average, and IF is the fairness concept impacting the worst-case utilitarian outcome of a lottery the most.

\begin{figure}[!ht]
    \centering
    \begin{subfigure}{0.45\textwidth}
        \includegraphics[width=\textwidth]{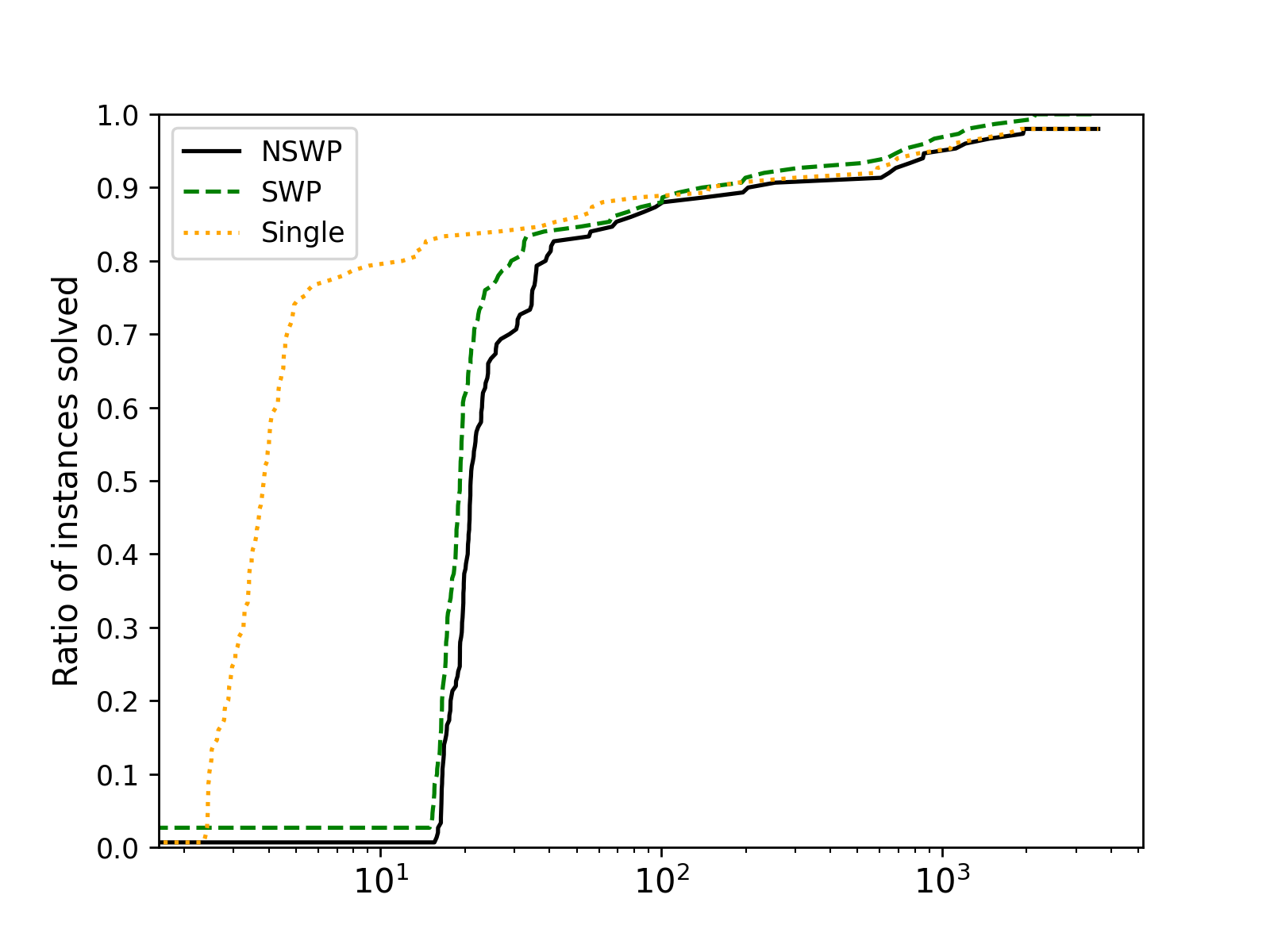}
        \caption{IF}
    \end{subfigure}
    \begin{subfigure}{0.45\textwidth}
        \includegraphics[width=\textwidth]{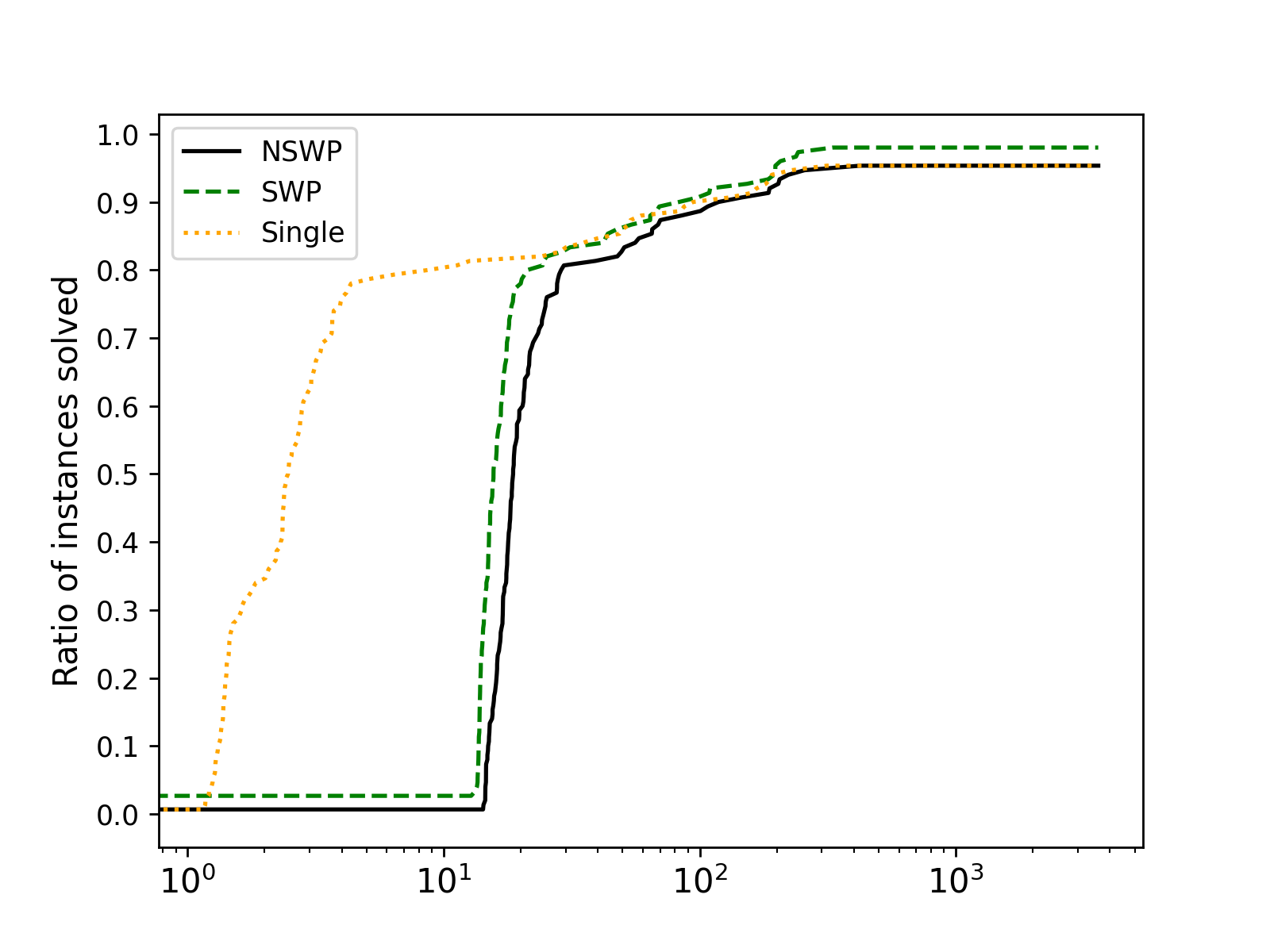}
        \caption{Rawls}
    \end{subfigure}
    
    \begin{subfigure}{0.45\textwidth}
        \includegraphics[width=\textwidth]{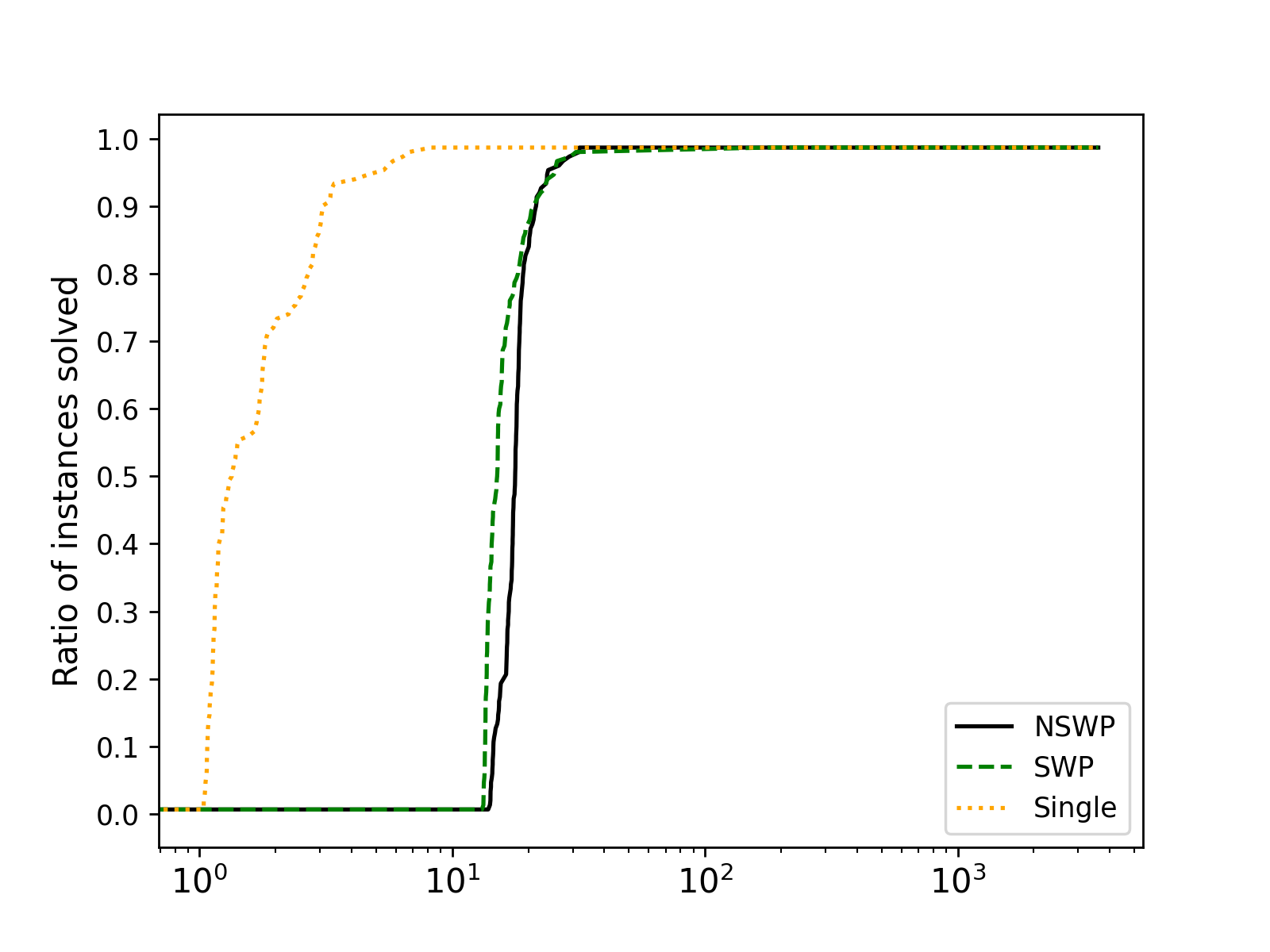}
        \caption{Aristotle}
    \end{subfigure}
    \begin{subfigure}{0.45\textwidth}
        \includegraphics[width=\textwidth]{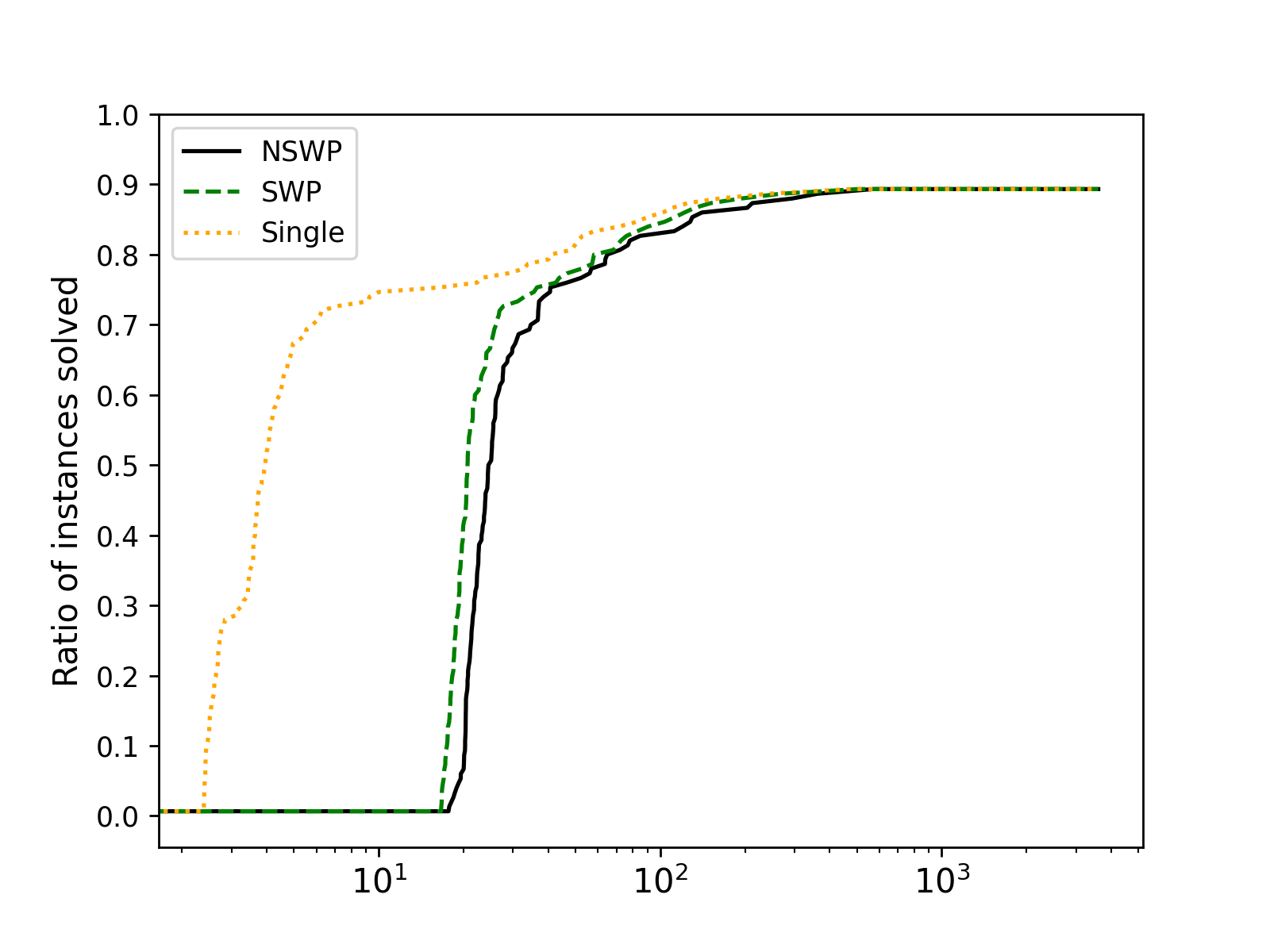}
        \caption{Nash}
    \end{subfigure}
    \caption{Ratio of instances solved for each model over time}
    \label{fig:times_scheme}
\end{figure}

\begin{figure}[!ht]
    \centering
    \begin{subfigure}{0.45\textwidth}
        \centering
        \includegraphics[width=\textwidth]{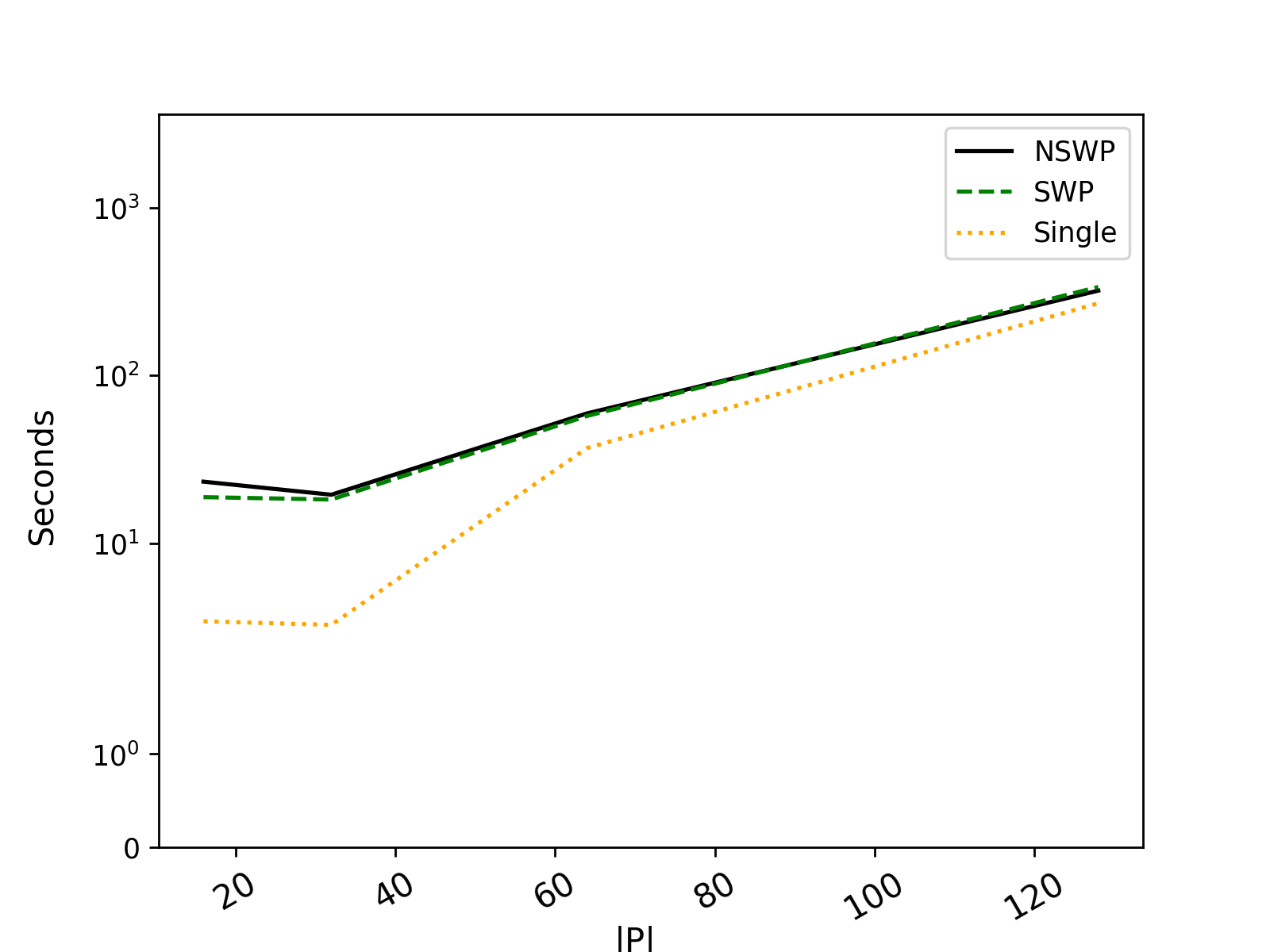}
        \caption{IF}
    \end{subfigure}
    \begin{subfigure}{0.45\textwidth}
        \centering
        \includegraphics[width=\textwidth]{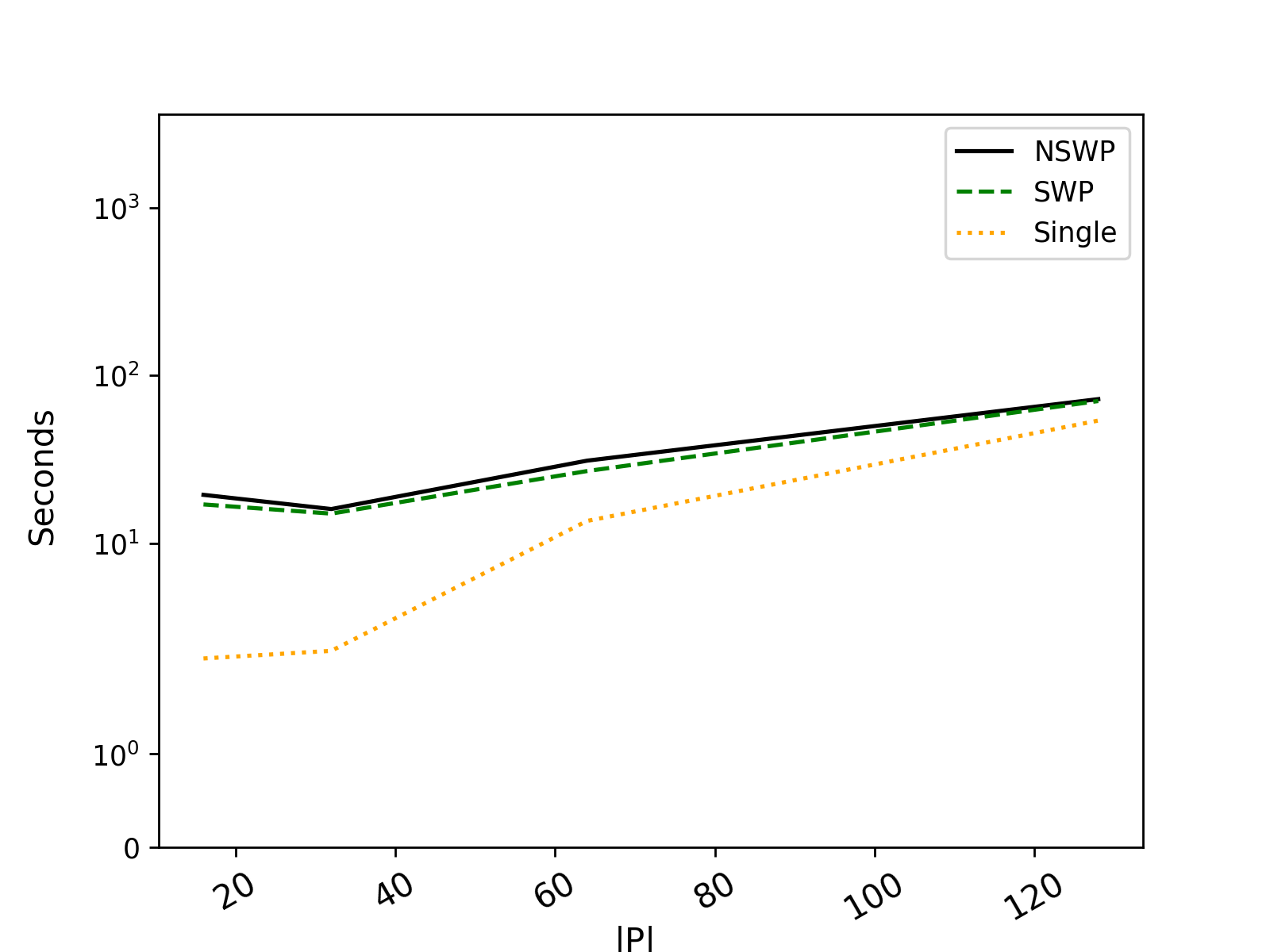}
        \caption{Rawls}
    \end{subfigure}
    
    \begin{subfigure}{0.45\textwidth}
        \centering
        \includegraphics[width=\textwidth]{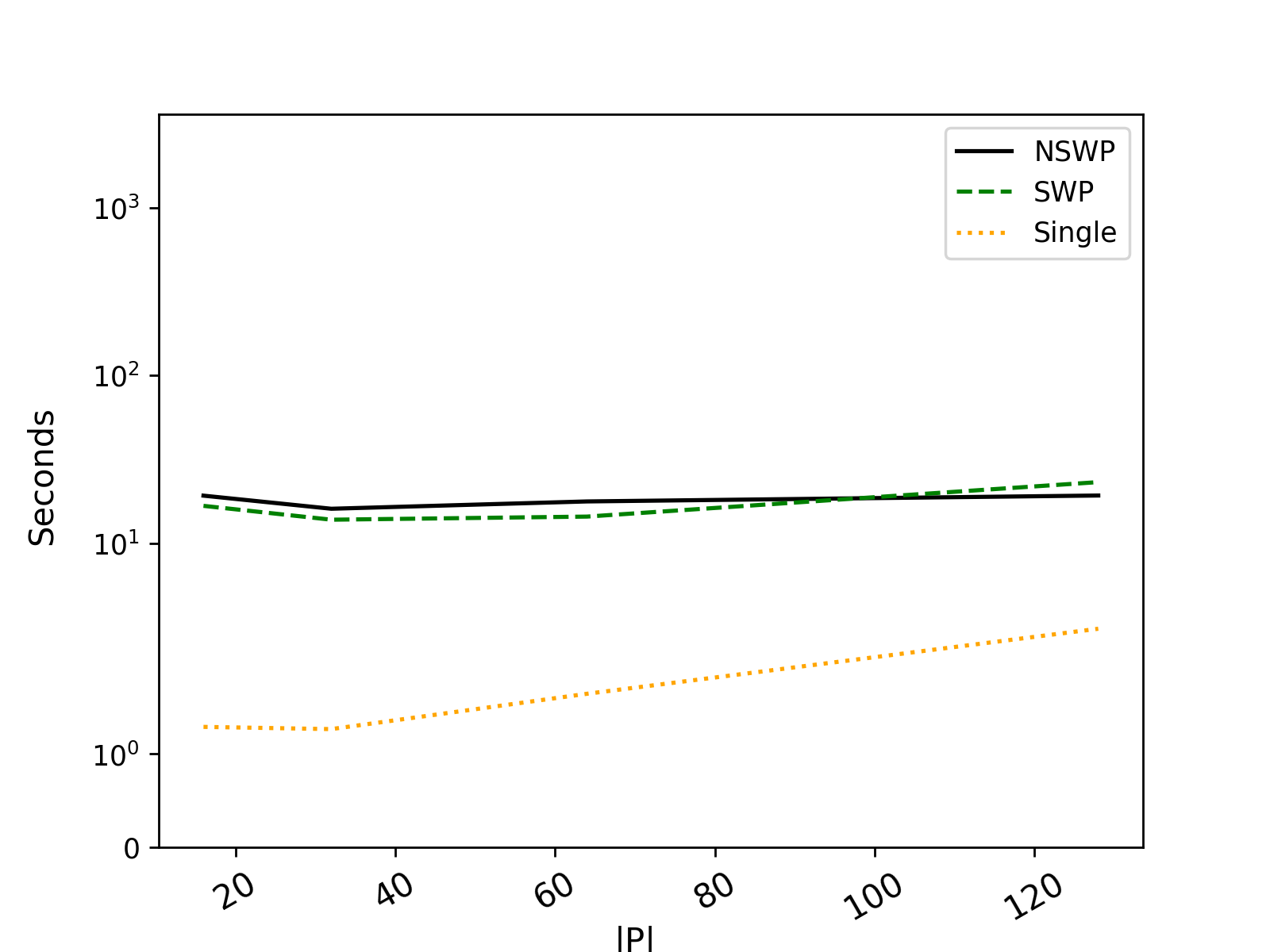}
        \caption{Aristotle}
    \end{subfigure}
    \begin{subfigure}{0.45\textwidth}
        \centering
        \includegraphics[width=0.98\textwidth]{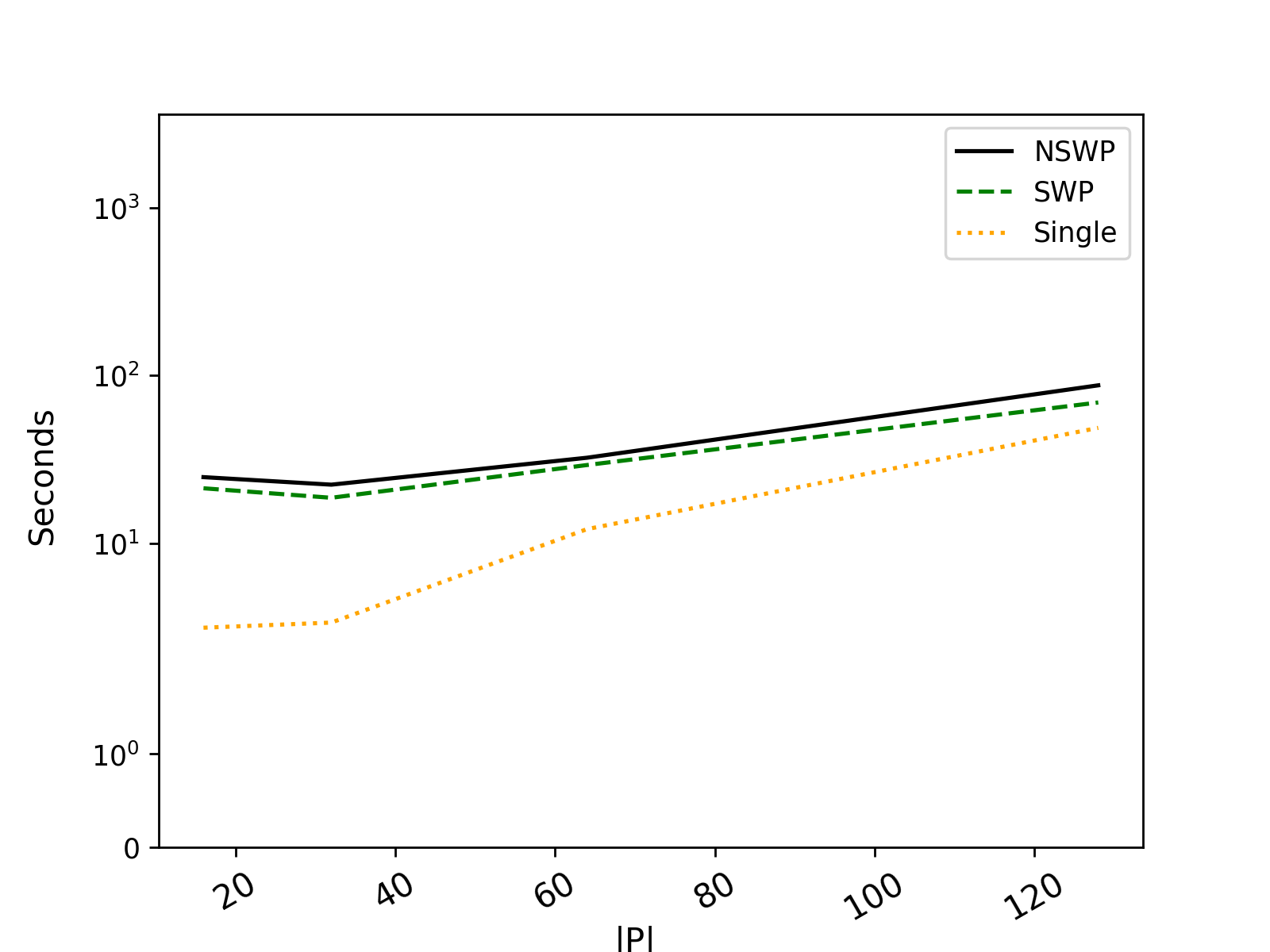}
        \caption{Nash}
    \end{subfigure}
    \caption{Time to solve the models for each scheme}
    \label{fig:times_size}
\end{figure}

\noindent\textbf{RQ4: How scalable are our schemes?} We compute the running time to solve Single, Linear and NSWP using our column generation approach. In Figure~\ref{fig:times_scheme}, we report performance profiles 
and in Figure~\ref{fig:times_size}, we  provide the running times averaged over the various sizes of $P$ for each scheme.

In all the fairness schemes, the curves for SWP and NSWP (see Figure~\ref{fig:times_scheme}) are very close and follow similar trajectories.  Figure~\ref{fig:times_size} underscores the same point: the time taken by SWP and NSWP is distributed very similarly. Both approaches are not empirically costly to run when compared to the single objective setting and given that this is a planning problem.

\noindent\textbf{RQ5: How many feasible exchange plans are part of the support of the optimal distribution?} We are interested in determining the number of feasible exchange plans in the support of the lotteries for our fairness schemes. This will complement our understanding of their running times. 

\begin{figure}[!ht]
    \centering
    \begin{subfigure}{0.45\textwidth}
        \centering
        \includegraphics[width=\textwidth]{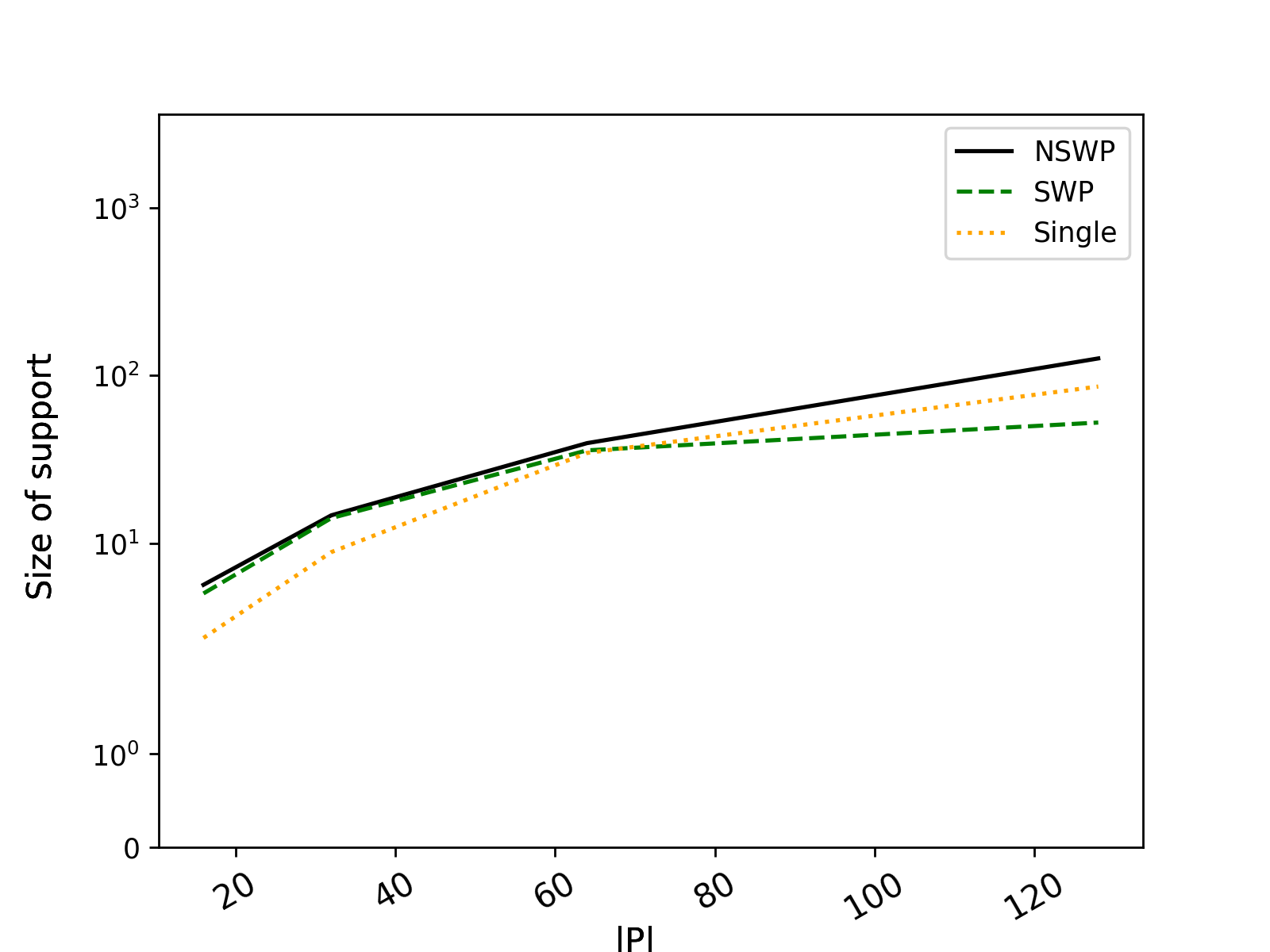}
        \caption{IF}
    \end{subfigure}
    \begin{subfigure}{0.45\textwidth}
        \centering
        \includegraphics[width=\textwidth]{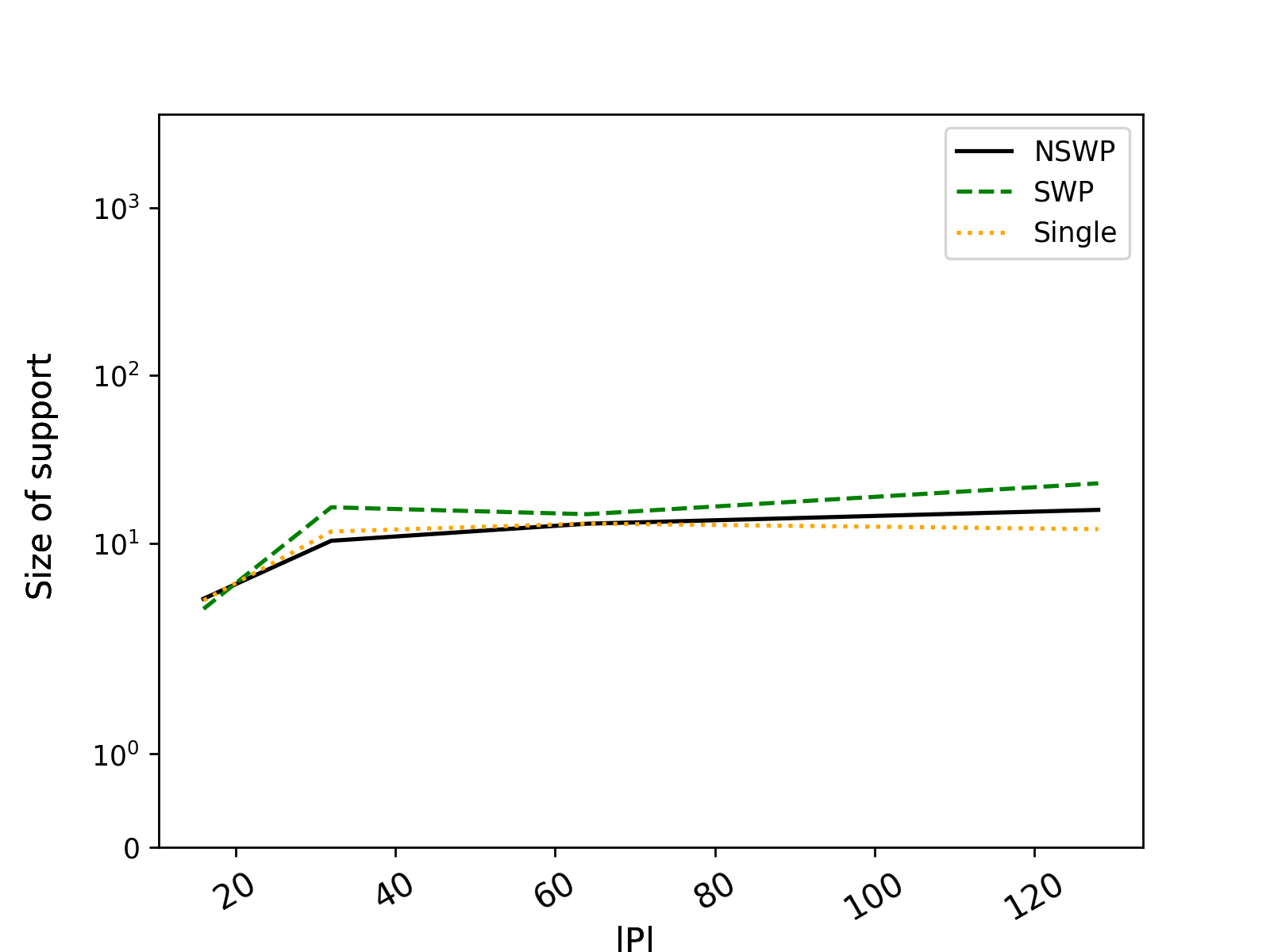}
        \caption{Rawls}
    \end{subfigure}
    
    \begin{subfigure}{0.45\textwidth}
        \centering
        \includegraphics[width=\textwidth]{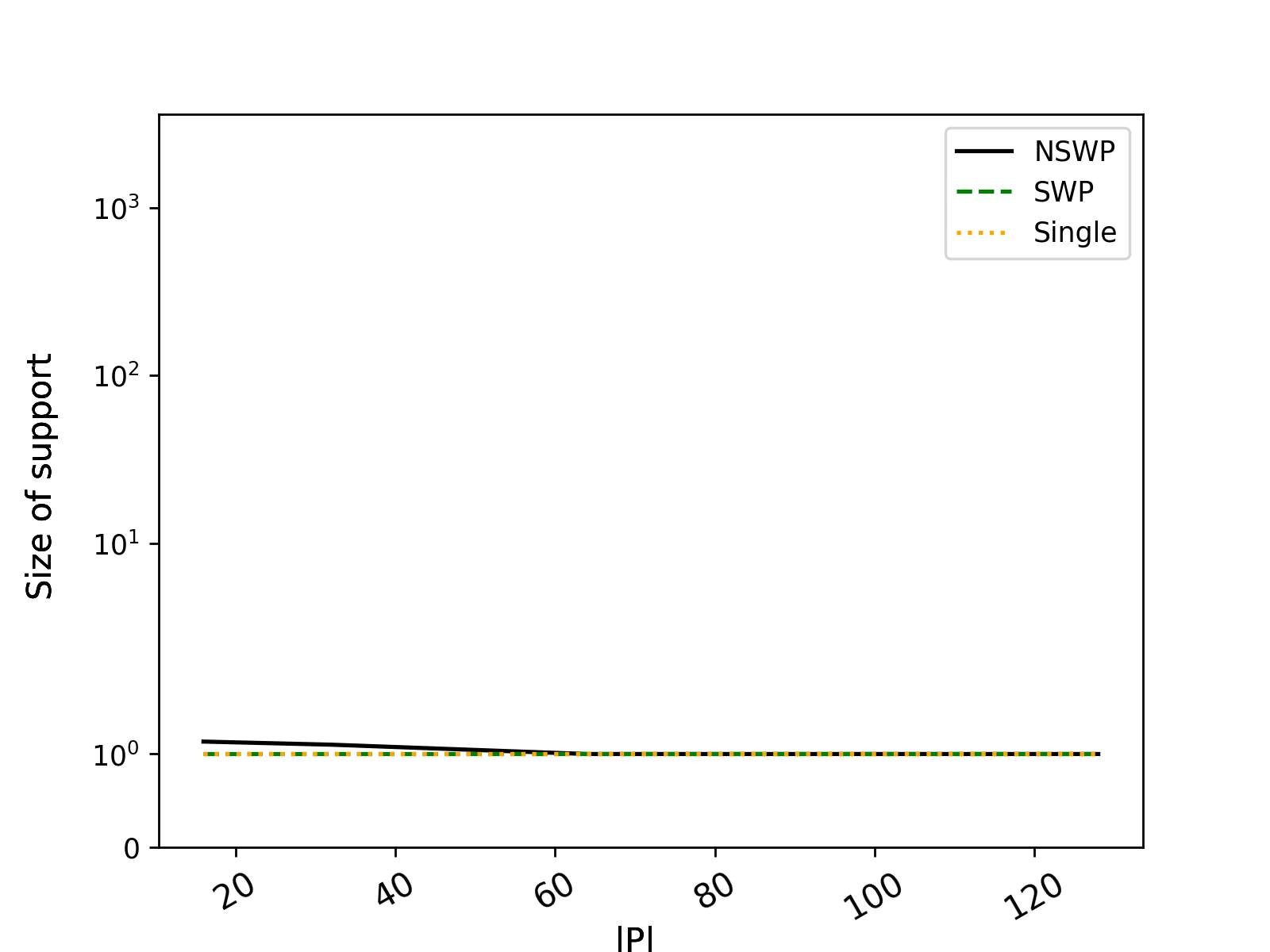}
        \caption{Aristotle}
    \end{subfigure}
    \begin{subfigure}{0.45\textwidth}
        \centering
        \includegraphics[width=0.98\textwidth]{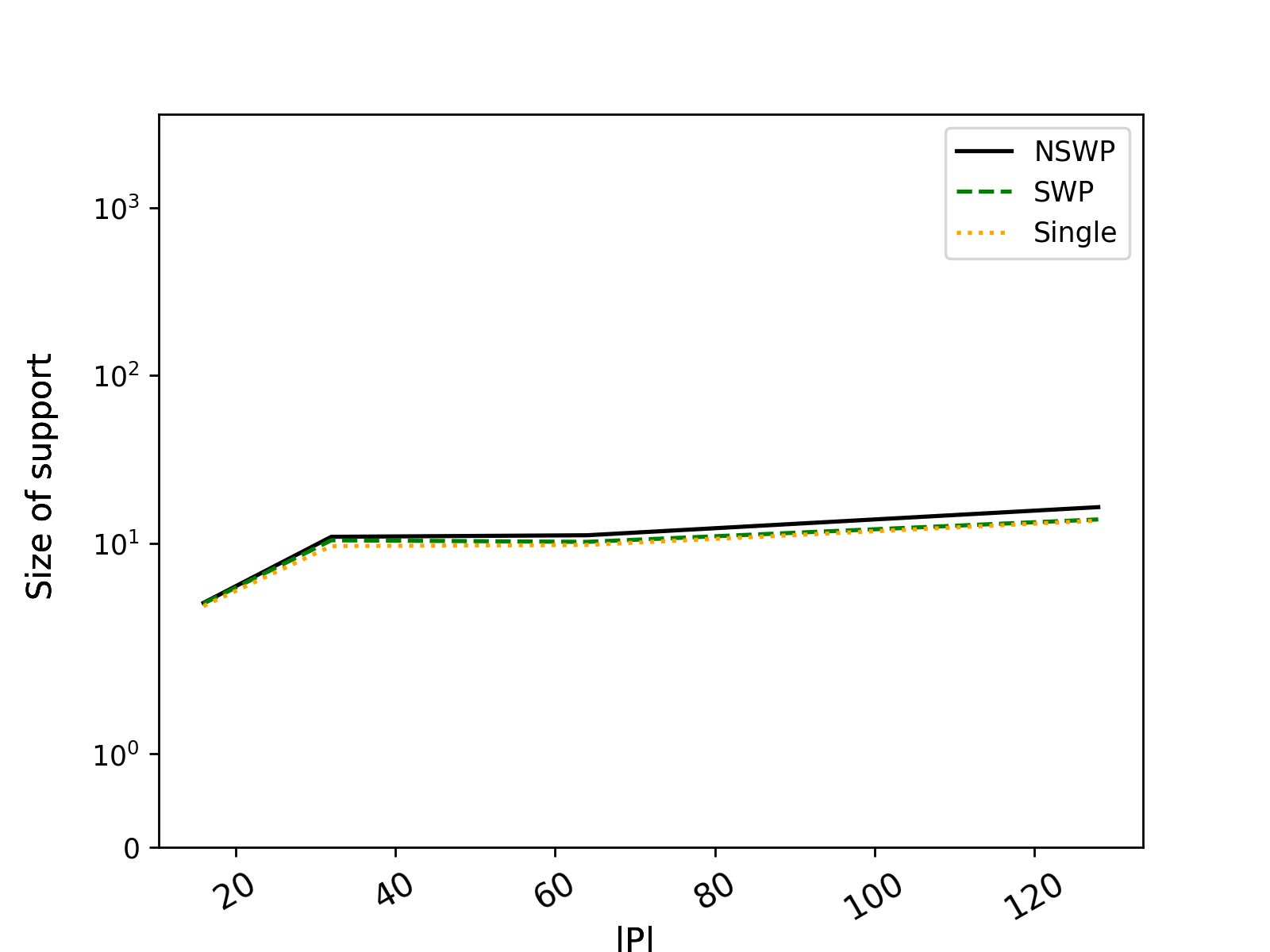}
        \caption{Nash}
    \end{subfigure}
    \caption{Size of the support for each scheme}
    \label{fig:support}
\end{figure}

We do note that the enumeration of exchange plans (columns) is usually stopped when there is a marginal increase to the objective (using our tolerance parameter). Therefore, we only have an approximation to the optimal distribution, but this is to be expected in any case because we are dealing with numerical methods. In Figure~\ref{fig:support}, we provide the average number of exchange plans in the support over all graph sizes for the three optimization approaches. We can see that the support size  is similar for all schemes (with a caveat for IF). The support size for Single is lower than NSWP and SWP for IF and Rawls. This intuitively makes sense as the multi-objective setting needs to include solutions that are in the support of the optimal distribution with respect to each objective. Schemes using Aristotle's objective have support size often equal to 1, because it is often possible to have a single solution that includes the maximum number of transplants and the maximum number of hard-to-match patients.

\section{Conclusions and future work}\label{sec:conclusions}
In this work, we have shown that fairness in KEPs inherently involves a tradeoff between maximizing the number of transplants and enhancing the chance of being selected in an exchange plan for patients that would otherwise be overlooked. We attempted to bridge these two dichotomies using the SWP and NSWP. We implemented this idea using various fairness schemes and provided a column generation for their efficient use. Our empirical results showed that the schemes based on the SWP and the NSWP fairness schemes have a low POF  and the worst-case utilitarian outcome is higher than schemes based only on fairness concepts. This means that their practical application is likely to keep the utilitarian objective close to its maximum value while increasing the patients' chances of being transplanted. Please see Appendix~\ref{sec:appendix_remarks} for further remarks concerning the overall analysis of the experiments provided in this article.

Our work contributes with a catalog and a systematization of fairness concepts for KEPs. It also provides a novel combination of fairness schemes with lottery policies for KEPs, with the potential to impact the way exchange plans are chosen in practice. Interestingly, most of the fairness metrics used in this paper do not consider clinical features of the patients, except for the Aristotle's equity principle. The latter uses the patients' PRA which is a value broadly used in practice and thus ``accessible'' to practitioners. The remaining fairness concepts rely on the mathematical topology of the KEP graphs. In future work, it would be of interest to go beyond this and consider individual features of the patients such as time in dialysis and health status, to name a few. Moreover, future research must account for the fact that KEPs are not static and run over time. This increases the flexibility of applying fairness concepts that are currently myopic. The generalization of fairness schemes to dynamic KEPs will likely pose computational difficulties since it would require the consideration of future feasible exchange plans which  make the problem highly combinatorial. Lastly, the improvement of individual fairness is an important line of study. For instance, the restriction of this fairness scheme to maximal exchange plans would likely address the issue surveyed in our experiments (see Appendix~\ref{sec:appendix_bilevel}).

\section*{Acknowledgements}
This work was funded by the FRQ-IVADO Research Chair in Data Science for Combinatorial Game Theory, NSERC grant 2019-04557 and 2021-04378, Canada CIFAR AI Chair, and Google scholar award. This research was enabled in part by support provided by  Calcul Qu\'ebec (\url{www.calculquebec.ca}) and Compute Canada (\url{www.computecanada.ca}).






\bibliographystyle{informs2014}
\bibliography{mybibfilearxiv}

\input{SupplementaryMaterial_IJOC}

\end{document}

%% file: SupplementaryMaterial_IJOC.tex
\newpage

\renewcommand{\thesection}{\Roman{section}} 
\renewcommand{\thesubsection}{\thesection.\Roman{subsection}}
\setcounter{section}{0}
\counterwithin{figure}{section}
\counterwithin{table}{section}

\begin{center}
{\Large Appendix}
\end{center}

\section{Conic programming definitions}
\label{sec:appendix_def}

\begin{deff}
    The second-order cone $\mathcal{Q}^n$ is given by
    \begin{equation*}
        \mathcal{Q}^n = \{x \in \mathbb{R}^n \mid x_1 \geq \sqrt{x_2^2 + \dots + x_n^2} \}
    \end{equation*}
\end{deff}
\begin{deff}
    The rotated second-order cone $\mathcal{Q}^n_{r}$ is given by
    \begin{equation*}
        \mathcal{Q}^n_{r} = \{x \in \mathbb{R}^n \mid 2 x_1 x_2 \geq x_3^2 + \dots + x_n^2, x_1, x_2 \geq 0 \}
    \end{equation*}
\end{deff}

\begin{deff}
The exponential cone $K_{\text{exp}}$ is defined by
\begin{equation*}
    K_{\text{exp}} = \left\{ x \in \mathbb{R}^3 \mid x_1 \geq x_2 e^{\frac{x_3}{x_2}}, x_2 > 0 \right\} \cup \left\{ x \in \mathbb{R}^3 \mid x_1 \geq 0, x_2 = 0, x_3 \leq 0  \right\}
\end{equation*}
and its dual $(K_{\text{exp}})^*$ is given by
\begin{equation*}
    (K_{\text{exp}})^* = cl \left\{ y \in \mathbb{R}^3 \mid y_1 \geq -y_3 e^{\frac{y_2}{y_3} - 1}, y_1 > 0, y_3 < 0 \right\}
\end{equation*}
\end{deff}

\section{Other formulations of the NSWP}
In this section of the appendix, we provide the missing dual formulations to the other primal problems introduced in the main body of the article. 

\subsection{NSWP with Rawlsian objective}

\begin{align*}
    \label{eq:dual_NSWP_mp}
    \tag{$\text{D}_{\text{Rawls}}$}
    \min -\alpha_0 + \alpha_1 d_1 + \alpha_2 d_2 \\
    \text{s.t.} \quad
    \begin{pmatrix}
        \alpha_1 \\
        \alpha_2 \\
        1
    \end{pmatrix} + u &= 0 \\
    \beta_v + \gamma_v &= 0 &\forall v \in P'\\
    -\alpha_2 - \eta &= 0\\
    \alpha_0 - \sum_{v \in V(S) \setminus N} (\alpha_1 + \beta_v) + \lambda_S &= 0 &\forall S \in \mathcal{F}_G \\
    u &\in \mathcal{Q}_{r}^3 \\
    \gamma_v &\geq 0 &\forall v \in P' \\
    \alpha \in \mathbb{R}^3, \beta \in \mathbb{R}^{|P'|}, \eta \in \mathbb{R}, \lambda &\geq 0.
\end{align*}

\subsection{NSWP with Aristotelian objective}

\begin{align*}
    \label{eq:dual_NSWP_aristotle}
    \tag{$\text{D}_{\text{Aristotle}}$}
    \min -\alpha_0 + \alpha_1 d_1 + \alpha_2 d_2 \\
    \text{s.t.} \quad
    \begin{pmatrix}
        \alpha_1 \\
        \alpha_2 \\
        1
    \end{pmatrix} + u &= 0 \\
    \alpha_0 - \sum_{v \in V(S) \setminus N} (\alpha_1 + \alpha_2) + \lambda_S &= 0 &\forall S \in \mathcal{F}_G \\
    u &\in \mathcal{Q}_{r}^3 \\
    \alpha \in \mathbb{R}^3, \lambda &\geq 0.
\end{align*}

\subsection{NSWP with Nash principle of fairness}

\begin{align*}
    \label{eq:dual_NSWP_nash}
    \tag{$\text{D}_{\text{Nash}}$}
    \min -\alpha_0 + \alpha_1 d_1 + \alpha_2 d_2 + \sum_{v \in P'} (w_v)_2 \\
    \text{s.t.} \quad
    \begin{pmatrix}
        \alpha_1 \\
        \alpha_2 \\
        1
    \end{pmatrix} + u &= 0 \\
    \begin{pmatrix}
        \beta_v \\
        \eta
    \end{pmatrix} + 
    \begin{pmatrix}
        (w_v)_1 \\
        (w_v)_3
    \end{pmatrix} &= 0 &\forall v \in P'\\
    -\alpha_2 - \eta &= 0\\
    \alpha_0 - \sum_{v \in V(S) \setminus N} (\alpha_1 + \beta_v) + \lambda_S &= 0 &\forall S \in \mathcal{F}_G \\
    u &\in \mathcal{Q}_{r}^3 \\
    w_v &\in (K_{\text{exp}})^* &\forall v \in P' \\
    \alpha \in \mathbb{R}^4, \beta \in \mathbb{R}^{|P'|}, \eta \in \mathbb{R}, \lambda &\geq 0.
\end{align*}

\section{Reference point}
\label{sec:appendix_reference}
Next, we detail the computation of the reference  point $(d_1,d_2)$, essential for the formulation of the NSWP schemes. We begin by giving a definition of the ideal and reference points in MOO.

\begin{deff}
Given a vector-valued function $f: X \to \mathbb{R}^k$, we define the ideal $i$ to be 
\begin{equation*}
    i =
    \begin{pmatrix}
        \sup_{x \in X} f_1(x) \\
        \vdots \\
        \sup_{x \in X} f_k(x).
    \end{pmatrix}
\end{equation*}
\label{def:ideal}
\end{deff}

\begin{deff}
Given a vector-valued function $f: X \to \mathbb{R}^k$, we define the nadir $d$ to be 
\begin{equation*}
    n =
    \begin{pmatrix}
        \inf_{x \in X} f_1(x) \\
        \vdots \\
        \inf_{x \in X} f_k(x).
    \end{pmatrix}
\end{equation*}
\label{def:nadir}
\end{deff}

In order to obtain the reference point for our NSWP models with respect to IF, Rawls and Aristotle, we partially build the models without including the terms $d_1, d_2$ and we use only $y_2$ as an objective. We can use Algorithm~1 to generate exchange plans that will lead to the best $f_2$ score possible. We then set this value to be $i_2$. We add the temporary constraint $y_2 = i_2$, replace the objective by $y_1$ and solve the mathematical program, obtaining $d_1$. Afterwards, we  replace the temporary constraint by $y_1 = i_1$, where $i_1$ is the maximum number of transplants that can be performed. Solving the model with objective $y_2$, we obtain $d_2$. Thus, the nadir point is given by $(d_1, d_2)$ (in the schemes that use it, \ie, IF, Rawls, and Aristotle), while the ideal is $(i_1, i_2)$. 

Nash's principle of fairness is slightly different. Since it sometimes occurs that some vertices are never included by a solution that maximizes the number of transplants, the value $d_2$ would be negative infinity. This is a problem and we must therefore rely on an ad hoc notion. This notion gives the worst-case scenario that the best value of $f_2$ can take over all graphs. In other words, if we have a complete bipartite graph $G = (V, A)$, where $V = P \cup N$ with $N = \{1\}$ and $P = \{2, \dots, n+1\}$, we can see that the only possible exchanges involve selecting vertex $1$ and one from $P$. In this case, the optimal solution selects the vertices from $P$ with probability $\frac{1}{n} = \frac{1}{\lvert P \rvert}$. Using this fact, we can set $d_2$ to be equal to $-\lvert P \rvert \log(\lvert P \rvert)$.

We can add the values of $d_1$ and $d_2$ to the right constraints and finalize the construction of our mathematical programming models. The exchange plans that were generated during this whole process can be used as a warm start to our set of solutions for Algorithm~1 when maximizing the NSWP models with objective $r$.

\section{Distance to the ideal vector in objective space}
\label{sec:appendix_dist}
An interesting property of solutions to the different approaches proposed in this article is to observe other measures of the quality of the selected solutions. Of course, these are measures that the method did not explicitly optimize and are intended to provide the reader with a more general understanding of what kind of solutions are returned by the SWP and the NSWP in practice.

For this purpose, we set out to answer the following research question:
\noindent\textbf{RQ6: How close is the optimal solution to the ideal vector in the objective space? How far is the optimal solution from the reference vector in the objective space?} To answer this question, we compute the relative distance from the ideal and reference points. It consists of the Euclidean distance from the ideal vector (resp., the reference vector) scaled by the distance between the ideal and reference vectors. These measures are averaged over all instances in order to better understand the relationship between distance to the ideal and reference vectors and how the two objective are balanced in the various approaches. See Tables~\ref{tab:ideal_dist} and~\ref{tab:nadir_dist}.

\begin{table}[!ht]
    \centering
    \begin{tabular}{c|c|c|c|c}
        \toprule
        & \multicolumn{4}{c}{Fairness concept} \\
    Scheme    & IF & Rawls & Aristotle & Nash\\
        \midrule
        SWP & $0.246 \pm 0.333$ & $0.259 \pm 0.368$ & $0.194 \pm 0.397$ & $0.002 \pm 0.006$\\
        NSWP & $0.256 \pm 0.332$ & $0.283 \pm 0.365$ & $0.194 \pm 0.397$ & $0.002 \pm 0.006$  \\
        \bottomrule
    \end{tabular}
    \caption{Relative distance from ideal vector for each scheme in the NSWP (all graph sizes)}
    \label{tab:ideal_dist}
\end{table}
\begin{table}[!ht]
    \centering
    \begin{tabular}{c|c|c|c|c}
        \toprule
        & \multicolumn{4}{c}{Scheme} \\
        & IF & Rawls & Aristotle & Nash\\
        \midrule
        SWP & $0.300 \pm 0.368$ & $0.213 \pm 0.335$ & $0.000 \pm 0.000$ & $0.998 \pm 0.005$\\
        NSWP & $0.286 \pm 0.348$ & $0.159 \pm 0.246$ & $0.000 \pm 0.000$ & $0.998 \pm 0.005$  \\
        \bottomrule
    \end{tabular}
    \caption{Relative distance from the reference vector for each scheme in the NSWP (all graph sizes)}
    \label{tab:nadir_dist}
\end{table}

SWP is both farther from the reference point and closer to the ideal than NSWP on average. These are proxy measures, however, as the NSWP is optimzing the area using the reference point and the solution as opposite corners. A smaller distance to the ideal can simply be the product of having close to the maximum number of transplants but performing poorly in terms of the fairness objective. The main takeaway is that the NSWP can sacrifice closeness to the ideal during its balancing of the objectives. This highlights the fact that when evaluating performance in a multi-objective sense, one needs to be careful about which criteria they care about. The distances to the ideal or from the reference point can certainly imply a balance of the objectives, but this cannot be guaranteed and we might in fact recover a solution that is not even on the Pareto frontier.

\section{Comparisons of various fairness schemes}
\label{sec:appendix_pou}
We proceed to a research question evaluating the performance of the NSWP in terms of fairness objectives. \textbf{RQ7}: How does the NSWP fare with respect to each fairness objective? Again, we consider IF, Rawlsian justice, Aristotle and Nash’s fairness principles.
\begin{figure}[!ht]
    \centering
    \begin{subfigure}{0.49\textwidth}
        \centering
        \includegraphics[width=\textwidth]{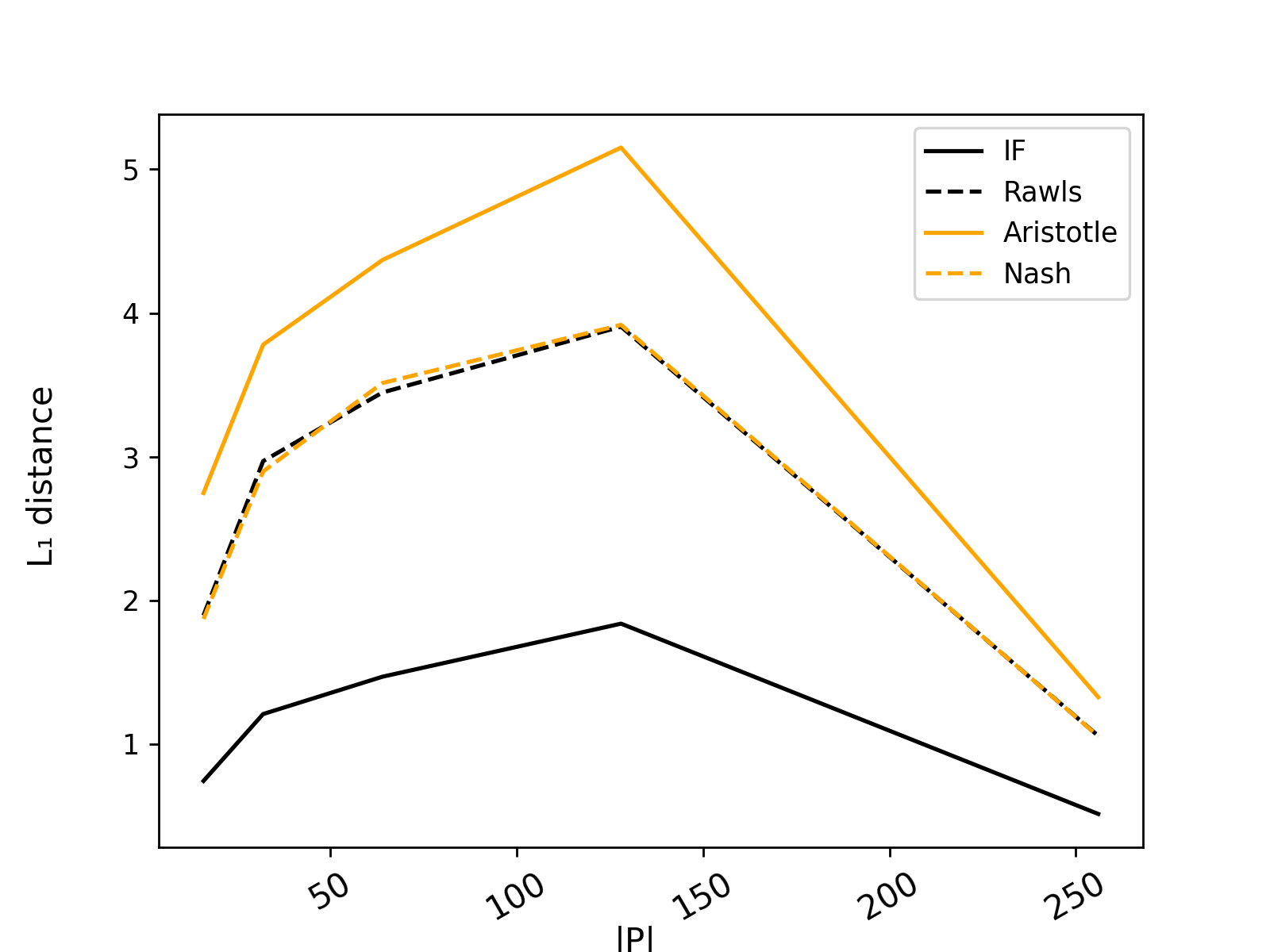}
        \caption{IF}
    \end{subfigure}
    \begin{subfigure}{0.49\textwidth}
        \centering
        \includegraphics[width=\textwidth]{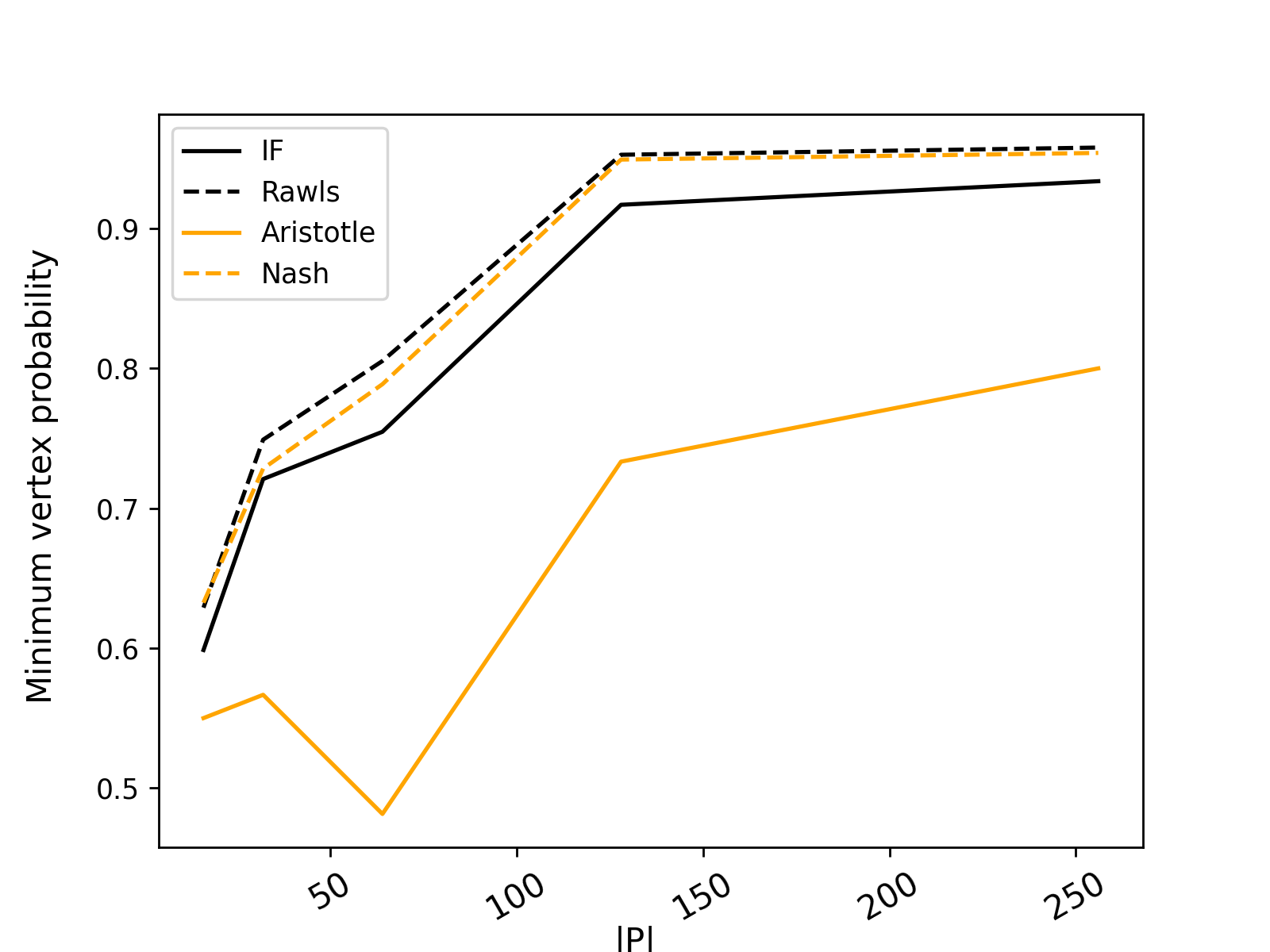}
        \caption{Rawls}
    \end{subfigure}
    
    \begin{subfigure}{0.49\textwidth}
        \centering
        \includegraphics[width=\textwidth]{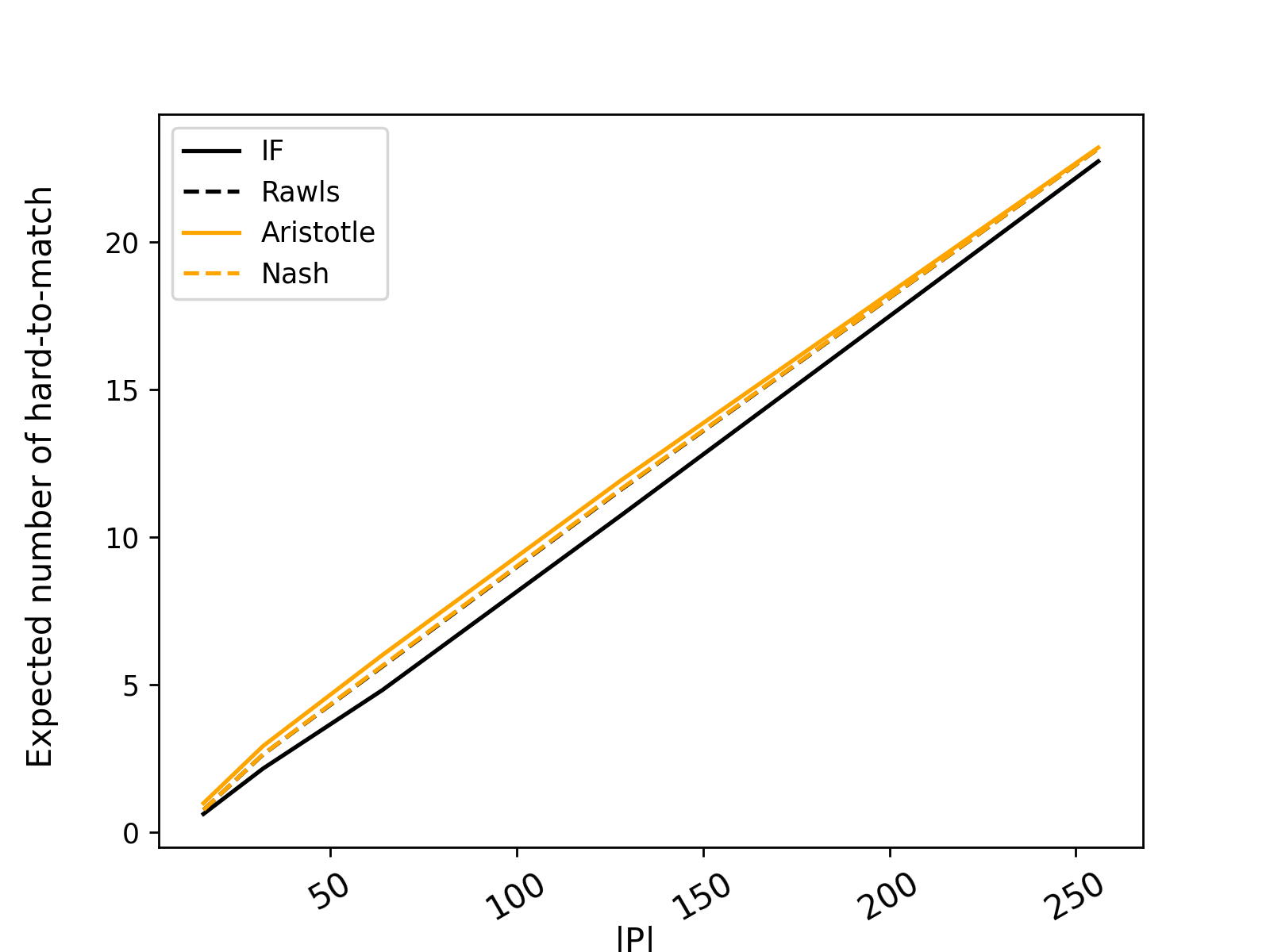}
        \caption{Aristotle}
    \end{subfigure}
    \begin{subfigure}{0.49\textwidth}
        \centering
        \includegraphics[width=0.98\textwidth]{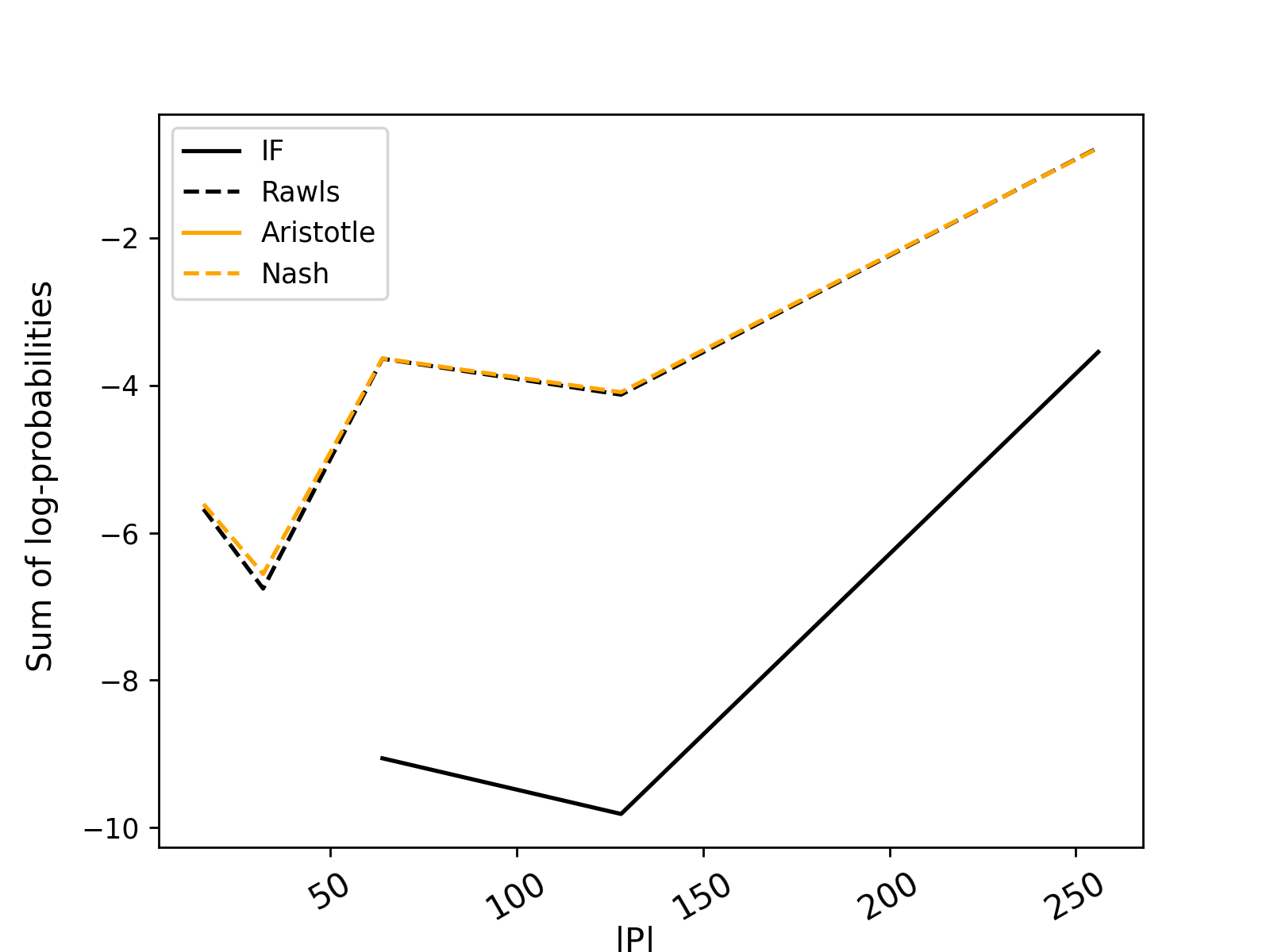}
        \caption{Nash}
    \end{subfigure}
    \caption{Comparison of the NSWP fairness schemes with respect to the values of each fairness objective}
    \label{fig:POU_comparison}
\end{figure}

By inspecting Figure~\ref{fig:POU_comparison}, we can observe that the (NSWP) Rawlsian and Nash fairness schemes perform very similarly. This is in sharp contrast to IF and Aristotle's fairness principle as they usually provide very different results. Aristotle's principle is noteworthy as it performs the worst on all objectives, except when it corresponds to the one being optimized.

\section{Issue related to IF}
\label{sec:appendix_bilevel}
The fairness scheme IF is problematic when we consider sub-optimal exchange plans in terms of the total number of transplants. This is because IF tries to minimize the distance to the average selection probability over all vertices in $P'$. However, if no vertex is selected, all probabilities are 0 and the ``best" solution is to select the empty exchange plan. When this is mixed with a utilitarian objective in a lottery setting over a graph with a single cycle, half of the time, we select the empty exchange plan, while we select the optimal exchange plan for the utilitarian objective the other half of the time (\ie, the single cycle). This leads to very poor results in practice and is due to the fact that IF is poorly designed for the inclusion of suboptimal exchange plans. We thus propose an interesting fix to this issue that involves bilevel programming. 

\begin{deff}
A feasible exchange plan $F \in \mathcal{F}_G$ is maximal if it is not contained in some strictly larger exchange plan $F'$. In other words, there does not exist an $F' \in \mathcal{F}_G$ such that $F$ is a subgraph of $F'$.
\end{deff}
Using this definition, we get rid of the optimality of the empty exchange plan for the fairness objective. However, because we add the constraint for maximal exchange plans, we must ensure that the generated solutions satisfy it. This leads to a bilevel formulation for the subproblem in the solution generation process.
\begin{align*}
    \label{eq:NSWP_subprob_bilevel}
    \tag{SP}
    \min \alpha_0 &- \sum_{v \in V(F) \cap P} (\alpha_1 + \alpha_3 + \beta_v) \\
    \text{s.t.} \quad 
    f_1(F') &\leq f_1(F) \\
    F' &\in \arg\max \left\{ f_1(S) :  S \in \mathcal{F}_G \text{ s.t. $F$ is a subgraph of $S$} \right\} \\
    F &\in \mathcal{F}_G
\end{align*}

The first constraint in \eqref{eq:NSWP_subprob_bilevel} ensures that we do not consider exchange plans that are not maximal and it will get rid of the aforementioned issue related to IF. It is an interesting research direction to explore how efficient this method can be made in practice and whether other improvements can be made to it.

\section{Instances solved}
\label{sec:appendix_solved}
In this section, we present the list of instances that were considered in our experiments. For each instance type (\ie, size of $P$ and percentage of NDDs), we specify how many models among Single, SWP and NSWP were able to terminate. We report these values in Tables~\ref{tab:instances_solved_single}, \ref{tab:instances_solved_swp} and \ref{tab:instances_solved_nswp} respectively. We can observe that both SWP and NSWP solve the vast majority of instances for which we were able to at least complete Single. Optimzing only the utilitarian objective did not result in more instances being solved; the bottleneck was thus related to the topology of the KEP instances.

\begin{table}[!ht]
    \centering
    \begin{subtable}{0.45\textwidth}
    \centering
    \begin{tabular}{c|c|c|c|c}
        \toprule
         &  \multicolumn{4}{c}{Percentage of NDDs}  \\
         $\lvert P \rvert$ & 0\% & 5\% & 10\% & 15 \% \\
         \midrule
         16 & 10 & 0 & 10 & 10\\
         32 & 10 & 10 & 10 & 10 \\
         64 & 10 & 10 & 10 & 10 \\
         128 & 10 & 9 & 10 & 10 \\
         256 & 10 & 10 & 10 & 0 \\
         \bottomrule
    \end{tabular}
    \caption{IF}
    \end{subtable}
    \begin{subtable}{0.45\textwidth}
    \centering
    \begin{tabular}{c|c|c|c|c}
        \toprule
         &  \multicolumn{4}{c}{Percentage of NDDs}  \\
         $\lvert P \rvert$ & 0\% & 5\% & 10\% & 15 \% \\
         \midrule
         16 & 10 & 0 & 10 & 10\\
         32 & 10 & 10 & 10 & 10 \\
         64 & 10 & 10 & 10 & 10 \\
         128 & 10 & 9 & 10 & 10 \\
         256 & 10 & 10 & 10 & 0 \\
         \bottomrule
    \end{tabular}
    \caption{Rawls}
    \end{subtable}
    \begin{subtable}{0.45\textwidth}
    \centering
    \begin{tabular}{c|c|c|c|c}
        \toprule
         &  \multicolumn{4}{c}{Percentage of NDDs}  \\
         $\lvert P \rvert$ & 0\% & 5\% & 10\% & 15 \% \\
         \midrule
         16 & 10 & 0 & 10 & 10\\
         32 & 10 & 10 & 10 & 10 \\
         64 & 10 & 10 & 10 & 10 \\
         128 & 10 & 9 & 10 & 10 \\
         256 & 10 & 10 & 10 & 0 \\
         \bottomrule
    \end{tabular}
    \caption{Aristotle}
    \end{subtable}
    \begin{subtable}{0.45\textwidth}
    \centering
    \begin{tabular}{c|c|c|c|c}
        \toprule
         &  \multicolumn{4}{c}{Percentage of NDDs}  \\
         $\lvert P \rvert$ & 0\% & 5\% & 10\% & 15 \% \\
         \midrule
         16 & 10 & 0 & 10 & 10\\
         32 & 10 & 10 & 10 & 10 \\
         64 & 10 & 10 & 10 & 10 \\
         128 & 10 & 10 & 10 & 10 \\
         256 & 10 & 9 & 10 & 0 \\
         \bottomrule
    \end{tabular}
    \caption{Nash}
    \end{subtable}
    \caption{Number of instances solved for each combination of $\lvert P \rvert$ and percentage of NDDs (for each fairness scheme) for a single (fairness) objective}.
    \label{tab:instances_solved_single}
\end{table}

\begin{table}[!ht]
    \centering
    \begin{subtable}{0.45\textwidth}
    \centering
    \begin{tabular}{c|c|c|c|c}
        \toprule
         &  \multicolumn{4}{c}{Percentage of NDDs}  \\
         $\lvert P \rvert$ & 0\% & 5\% & 10\% & 15 \% \\
         \midrule
         16 & 9 & 0 & 9 & 10\\
         32 & 10 & 10 & 10 & 10 \\
         64 & 10 & 10 & 10 & 10 \\
         128 & 10 & 9 & 10 & 10 \\
         256 & 10 & 10 & 10 & 0 \\
         \bottomrule
    \end{tabular}
    \caption{IF}
    \end{subtable}
    \begin{subtable}{0.45\textwidth}
    \centering
    \begin{tabular}{c|c|c|c|c}
        \toprule
         &  \multicolumn{4}{c}{Percentage of NDDs}  \\
         $\lvert P \rvert$ & 0\% & 5\% & 10\% & 15 \% \\
         \midrule
         16 & 9 & 0 & 9 & 10\\
         32 & 10 & 10 & 10 & 10 \\
         64 & 10 & 10 & 10 & 10 \\
         128 & 10 & 9 & 10 & 10 \\
         256 & 10 & 10 & 10 & 0 \\
         \bottomrule
    \end{tabular}
    \caption{Rawls}
    \end{subtable}
    \begin{subtable}{0.45\textwidth}
    \centering
    \begin{tabular}{c|c|c|c|c}
        \toprule
         &  \multicolumn{4}{c}{Percentage of NDDs}  \\
         $\lvert P \rvert$ & 0\% & 5\% & 10\% & 15 \% \\
         \midrule
         16 & 9 & 0 & 9 & 10\\
         32 & 10 & 10 & 10 & 10 \\
         64 & 10 & 10 & 10 & 10 \\
         128 & 10 & 9 & 10 & 10 \\
         256 & 10 & 10 & 10 & 0 \\
         \bottomrule
    \end{tabular}
    \caption{Aristotle}
    \end{subtable}
    \begin{subtable}{0.45\textwidth}
    \centering
    \begin{tabular}{c|c|c|c|c}
        \toprule
         &  \multicolumn{4}{c}{Percentage of NDDs}  \\
         $\lvert P \rvert$ & 0\% & 5\% & 10\% & 15 \% \\
         \midrule
         16 & 9 & 0 & 9 & 10\\
         32 & 10 & 10 & 10 & 10 \\
         64 & 10 & 10 & 10 & 10 \\
         128 & 10 & 10 & 10 & 10 \\
         256 & 10 & 9 & 10 & 0 \\
         \bottomrule
    \end{tabular}
    \caption{Nash}
    \end{subtable}
    \caption{Number of instances solved for each combination of $\lvert P \rvert$ and percentage of NDDs (for each fairness scheme) by the SWP}
    \label{tab:instances_solved_swp}
\end{table}

\begin{table}[!ht]
    \centering
    \begin{subtable}{0.45\textwidth}
    \centering
    \begin{tabular}{c|c|c|c|c}
        \toprule
         &  \multicolumn{4}{c}{Percentage of NDDs}  \\
         $\lvert P \rvert$ & 0\% & 5\% & 10\% & 15 \% \\
         \midrule
         16 & 9 & 0 & 9 & 10\\
         32 & 10 & 10 & 10 & 10 \\
         64 & 10 & 10 & 10 & 10 \\
         128 & 10 & 9 & 10 & 10 \\
         256 & 10 & 10 & 10 & 0 \\
         \bottomrule
    \end{tabular}
    \caption{IF}
    \end{subtable}
    \begin{subtable}{0.45\textwidth}
    \centering
    \begin{tabular}{c|c|c|c|c}
        \toprule
         &  \multicolumn{4}{c}{Percentage of NDDs}  \\
         $\lvert P \rvert$ & 0\% & 5\% & 10\% & 15 \% \\
         \midrule
         16 & 9 & 0 & 9 & 10\\
         32 & 10 & 10 & 10 & 10 \\
         64 & 10 & 10 & 10 & 10 \\
         128 & 10 & 9 & 10 & 10 \\
         256 & 10 & 10 & 10 & 0 \\
         \bottomrule
    \end{tabular}
    \caption{Rawls}
    \end{subtable}
    \begin{subtable}{0.45\textwidth}
    \centering
    \begin{tabular}{c|c|c|c|c}
        \toprule
         &  \multicolumn{4}{c}{Percentage of NDDs}  \\
         $\lvert P \rvert$ & 0\% & 5\% & 10\% & 15 \% \\
         \midrule
         16 & 9 & 0 & 9 & 10\\
         32 & 10 & 10 & 10 & 10 \\
         64 & 10 & 10 & 10 & 10 \\
         128 & 10 & 9 & 10 & 10 \\
         256 & 10 & 10 & 10 & 0 \\
         \bottomrule
    \end{tabular}
    \caption{Aristotle}
    \end{subtable}
    \begin{subtable}{0.45\textwidth}
    \centering
    \begin{tabular}{c|c|c|c|c}
        \toprule
         &  \multicolumn{4}{c}{Percentage of NDDs}  \\
         $\lvert P \rvert$ & 0\% & 5\% & 10\% & 15 \% \\
         \midrule
         16 & 9 & 0 & 9 & 10\\
         32 & 10 & 10 & 10 & 10 \\
         64 & 10 & 10 & 10 & 10 \\
         128 & 10 & 10 & 10 & 10 \\
         256 & 10 & 9 & 10 & 0 \\
         \bottomrule
    \end{tabular}
    \caption{Nash}
    \end{subtable}
    \caption{Number of instances solved for each combination of $\lvert P \rvert$ and percentage of NDDs (for each fairness scheme) by the NSWP}
    \label{tab:instances_solved_nswp}
\end{table}

\section{Further remarks}
\label{sec:appendix_remarks}
After reviewing the various research questions and examining the results of the experiments, it might come in handy to explicitly describe individual instances and how the three methods fare with respect to the aforementioned statistics.

\begin{table}[!ht]
    \centering
    \resizebox{\linewidth}{!}{
    \begin{tabular}{l|c|c|c|c}
        \toprule
        & \multicolumn{4}{c}{Scheme} \\
        Statistic &  Single-objective (utility) & Single-objective (IF) & SWP & NSWP \\
        \midrule
        Solution values $(f_1,f_2)$ & $(16.000, -5.714) $ & $(0.000, 0.000)$ & $(12.000, -2.857)$ & $(10.500, -2.333)$\\
        Area & $0.000$ & $0.000$ & $34.286$ & $35.500$\\
        POF & $0.000$ & $1.000$ & $0.250$ & $0.344$ \\
        Distance to ideal & $5.714$ & $16.000$ & $4.916$ & $5.974$\\
        Distance to reference & $16.000$ & $5.714$ & $12.335$ & $11.030$\\
        \bottomrule
    \end{tabular}
    }
    \caption{Instance~32 with $(\lvert P \rvert,|N|)=(32,0)$  for IF}
    \label{tab:instance_32}
\end{table}

From Table~\ref{tab:instance_32}, we can make a key observation: the SWP approach can select a solution that is closer to the ideal vector than the NSWP. However, this solution does not Pareto dominate the solution from the NSWP. We can see that the total area given by the solution and the reference point for the NSWP is $35.5$ compared to $34.286$ for SWP. Hence, we can conclude that even if the solution to the SWP is close to the ideal vector, there is no guarantee that it provides a good balance of the two objectives. This point is further strengthened by the POF and $f_2$ values. The NSWP has a higher POF but its $f_2$ value is higher, again highlighting the balance of the two objectives.

From the results of our experiments, we have shown that the NSWP indeed balances utility and fairness while allowing for the expected number of transplants to be close to optimal. The key exception regards IF because the empty exchange plan can be selected as it is optimally fair according to the fairness objective. This is an issue that requires a novel column generation method approach for this particular case. Whether this can be effectively implemented in practice is beyond the scope of this paper. Our multi-objective approach does improve the fairness score of the distribution over feasible exchange plans when compared to the utilitarian objective case. However, it does so at the cost of additional computing time. From Figure~\ref{fig:POU}, we can see that for most instances, we do manage to obtain a result that improves the value of $f_2$ dramatically except for Aristotle's fairness principle and Nash's standard of comparison.

Aristotle's fairness principle did not yield interesting results in our experiments since we only optimized the number of patients with a PRA above 80\%. In the original definition of the POF, only the first solution returned by the solver was used to compute it. However, we have seen in practice that by including NDDs, we were able to generate a solution that both optimized the total number of transplants and includes the maximum number of hard-to-match patients possible in most cases. The fact that that there is an overlap between solutions maximizing the total number of transplants and those maximizing the number of hard-to-match patients is surprising. Nevertheless, we have to be careful as we are dealing with synthetic data. It is entirely possible that the way the graphs are generated is not sufficiently realistic. This suggests investigating the phenomenon more deeply.

